\documentclass[prc,aps,final,twocolumn,nofootinbib]{revtex4-1}
\usepackage{graphicx}
\usepackage{amsmath,amssymb,amsfonts}
\usepackage{bm}

\def\r{\mathbf{r}}
\def\p{\mathbf{p}}
\def\x{\xi}

\def\epsi{\mathcal{E}}
\def\cK{\mathcal{K}}

\def\Re{\mathop{\rm Re}}

\def\erf{\mathop{\rm erf}\nolimits}

\def\d{\hbox{d}}

\def\siml{\lesssim}
\def\be{\begin{equation}}
\def\ee{\end{equation}}
\def\bea{\begin{eqnarray}}
\def\eea{\end{eqnarray}}
\def\l{\label}
\def\sgn{\mathop{\rm sgn}}

\usepackage{colordvi}
\usepackage[normalem]{ulem}

\begin{document}

\title{Semiclassical catastrophe theory of simple bifurcations}

\author{A.\ G.\ Magner}
\email{Email: magner@kinr.kiev.ua}
\affiliation{\it Institute for Nuclear Research, 03680 Kyiv, Ukraine} 
\author{K.\ Arita}
\affiliation{\it Department of Physics, Nagoya Institute of 
Technology, Nagoya 466-8555, Japan} 
\date{\today}

\begin{abstract}
The Fedoriuk-Maslov catastrophe theory of caustics and turning points
is extended 
to solve the bifurcation problems by the improved
stationary phase method (ISPM).
The trace formulas for the radial power-law
(RPL) potentials are presented by the ISPM based on the second- and
third-order expansion of the classical action near the stationary
point.  A considerable enhancement of contributions of the two orbits
(pair of  
consisting of the parent and newborn 
orbits) at their bifurcation is shown. 
The ISPM trace formula is proposed for a simple bifurcation
scenario of Hamiltonian systems with continuous symmetries, 
where the contributions of the bifurcating parent orbits
vanish  
upon approaching the bifurcation point
due to the reduction of the end-point manifold.  This occurs since the
contribution of the parent orbits is included in  
the term
corresponding to the family of the newborn daughter orbits.
Taking this feature into account,
the ISPM level densities calculated for the RPL
potential model are shown to be in good agreement with the quantum
results at the bifurcations and asymptotically far from the
bifurcation points.
\end{abstract}
\maketitle

\section{INTRODUCTION}

Semiclassical periodic-orbit theory (POT) is a powerful tool for
the study of shell structures in the single-particle level density of
finite fermionic systems \cite{GUTZpr,GUTZbook90,BBann72,sclbook}.
This theory relates the oscillating level
density and shell-correction
energy to the sum over contributions of 
classical periodic orbits.
It thus gives the correspondence between
the fluctuation properties of the quantum dynamics and
the characteristics of the periodic motion embedded in the
classical dynamics.

Gutzwiller \cite{GUTZbook90} suggested the semiclassical evaluation of
the Green's function in the Feynman path integral representation to
derive the POT trace formula for the level density if the
single-particle Hamiltonian has no continuous symmetries other than
the time-translational invariance.  In this case, the energy $E$ is
the single integral of motion for a particle dynamics in the
mean-field potential.  For a given $E$, all the generic periodic
orbits (POs) are isolated, i.e., any variation of the initial
condition perpendicular to the PO will violate its periodicity. The
original version of the POT was extended for  
Hamiltonians with
continuous symmetries (extended Gutzwiller approach), in particular
for the rotational [$\mbox{O}(n)$] and oscillator-type [$\mbox{SU}(n)$]
symmetries \cite{SMepp76,SMODzp77,MAGosc78,CREAGHpra97}.
Berry and
Tabor \cite{BT76} derived the POT for integrable 
systems by applying the Poisson-summation method through the
semiclassical torus-quantization condition.  It is also helpful in the
case of a high classical degeneracy of POs.  Here the classical
degeneracy is defined by the number of independent parameters
$\mathcal{K}$ for a continuous family of 
classical periodic orbits
at a given energy of the particle.

Some applications of the POT to the deformations of nuclei and
metallic clusters by using 
phase-space variables were presented in
Refs. \cite{BT76,SMODzp77,CREAGHpra97,sclbook,migdalrev,MKApan16}.
The pronounced shell effects caused by the deformations have been
discussed.  Within the improved stationary-phase method
\cite{migdalrev,MFAMBpre01,MAFMptp02,MAFptp06,MFAMMSBptp99,KKMABps15,MKApan16}
(ISPM), the divergences and discontinuities of the
standard stationary-phase method 
(SSPM) \cite{BBann72,SMepp76,BT76,CREAGHpra97,sclbook}
near the symmetry-breaking
and bifurcation points were removed.

Bifurcations of the isolated POs in non-integrable Hamiltonian systems are
classified by the normal-form theory \cite{Ozoriobook,OH87},
based on a pioneering work of Meyer on the
differentiable symplectic mappings \cite{Meyer}.  The
change in the number of POs depends on the types of 
bifurcations:
from zero to two in the isochronous or saddle-node bifurcations,
from one to three in the period-doubling or pitchfork
bifurcations, and so on.  In the integrable 
systems, the classical phase space is entirely covered by the tori.
The classical orbits on the rational tori (in which the frequencies of
the independent motions are commensurable) form  
degenerate
families of POs.  In such systems, the bifurcations mostly occur at a
surface of the physical phase-space volume occupied by the classical
trajectories. For 
brevity, below, such a surface will be called  
the end point.  Note that even if the commensurability of the
frequencies is not fully satisfied at the end point, one has a PO
there because the mode with incommensurable frequencies has 
zero
amplitude.  In the action-angle representation $(I_n,\Theta_n)$ for
the $n$th mode, it corresponds to $I_n=0$, where the variation of
$\Theta_n$ will not generate a family.  Thus, the end-point POs have
smaller degeneracies than those inside the physical region.  With
varying  
potential parameter, a new PO family appears in the
transition from an unphysical to  
a physical region for a rational
torus where the frequencies of the end-point PO become commensurable.
It is considered as the bifurcation of the end-point PO generating a
new orbit family.  Such bifurcations take place repeatedly with
varying  
potential parameter and a new family of orbits with a
higher degeneracy and different frequency ratio is generated at each
bifurcation point.  Typical examples are the equatorial orbits with
$\cK=1$ in the spheroidal cavity where the generic family has $\cK=2$,
and the circular orbits with $\cK=2$ in a
spherical potential where the generic family appears with $\cK=3$.

Our ISPM is based on the catastrophe theory by Fedoriuk and Maslov for
solving problems with  
caustic and turning points in calculations
of the integrals by using the saddle-point
method \cite{FEDoryuk_pr,MASLOV,FEDoryuk_book1,MAFptp06}.  In the SSPM,
the catastrophe integrals are evaluated by an expansion of the action
integral in the exponent up to 
second-order terms and amplitude up
to zeroth-order  
terms near the stationary point, and the
integration limits are extended to  
an infinite interval.  The trace
formula based on the SSPM encounters a divergence or discontinuity at
the bifurcation point.  Such catastrophe problems are due to the zeros
(caustics) or infinities (turning points) of the second-order
derivatives of the 
action integral in this expansion.  Fedoriuk was the
first 
to prove \cite{FEDoryuk_pr,FEDoryuk_book1} the so-called
Maslov theorem \cite{MASLOV}: Each simple caustic or turning point
(having a finite nonzero third-order derivative of the action
integral) leads to a shift of the phase in the exponent of the catastrophe
integral by $-\pi/2$ along the classical trajectory.  This provides 
an extension of the one-dimensional WKB formula to higher dimensions.
Thus, in the asymptotic region far from the caustic and turning
points, one can use  
second-order expansions of the 
action integral
(and zeroth-order 
expansions of the amplitude) in the semiclassical Green's
function taking into account the shift of the phase according to the
Maslov theorem.  However, a proper foundation for an extension of the
Fedoriuk-Maslov catastrophe theory (FMCT) to the derivations of the
trace formulas near the PO bifurcations, where one has to treat the
semiclassical propagators near the catastrophe points, is still an
open question.

Uniform approximations based on the normal-form theories give
alternative ways of solving the bifurcation problems, see Refs.~
\cite{CREAGHpra97,sclbook,Ozoriobook,Ullmo,SIEjpa97,SSUNjpa97,
SCHOjpa98,ABjpa02,KAijmpe04,BO05,ABjpa08,KAprc12,MO15}.  The trace
formula valid near the bifurcation points is derived by calculating
the catastrophe integral in  
local uniform approximations.
It gives an indistinct combination of 
contributions of
bifurcating POs, and the
asymptotic regions (where each of the POs has a 
separate contribution to the
trace integral) are connected by a kind of the interpolation
through the bifurcations in the 
global uniform
approximations.  The ISPM
provides much simpler semiclassical formula in which one needs no
artificial interpolation procedures over the bifurcations.  Within the
simplest ISPM, the contributions of the bifurcating orbits are given
separately through the bifurcation 
and they are basically independent 
of the type of  
bifurcation.
This allows one to 
give analytic expressions for the Gaussian-averaged
level densities and the energy shell
corrections \cite{STRUTscm,FUHIscm1972,SMepp76,sclbook}.

In the present work we apply the FMCT  
to solve the bifurcation
problems 
that arise 
for some parameters of the mean-field potential
for the particle motion in the end-point phase space. 
For  
the simplest
but a rich exemplary case 
that nevertheless includes all of the
necessary points of the general behavior of the integrable systems, we 
will consider the spherical
radial power-law (RPL) potential $V\propto r^\alpha$
as a function of the radial power parameter $\alpha$, which
controls the surface diffuseness of the
system \cite{KAps16,KAijmpe04,KAprc12}.
 Some of these results are general for any integrable
and nonintegrable Hamiltonian systems, and can be applied to a more
realistic nuclear mean-field potential having the
deformation. 
The spherical RPL model has  
already been analyzed within the ISPM
in Ref.~\cite{MKApan16}. 
Good agreement between the semiclassical and quantum shell structures
was shown in the level-density and energy shell corrections 
for several values of the surface diffuseness parameter
including its symmetry-breaking and bifurcation values.
Quantum-classical correspondences in the deformed RPL models
with and without spin-orbit coupling were also studied and
various properties of the nuclear shape dynamics, 
such as the origins
of exotic deformations and the prolate-oblate asymmetries, have been
clarified \cite{KAprc12,AMprc14,KAps16,MKApan16,KAps17}.
Therefore, a proper study of the general aspects of the bifurcation problem
within the ISPM, even taking the simplest spherical RPL potential as an 
example
for which one can achieve 
much progress in analytical
derivations, is expected to be helpful.

This article is organized in the following way. In
Sec.~\ref{sec2:trace} we present a general semiclassical phase-space
trace formula for the level density as a typical catastrophe integral.
Section~\ref{sec3:cat3} shows the Fedoriuk-Maslov method for solving
the simple caustic- and turning-point singularity problems.
Section~\ref{sec4:bifurcations} is devoted to the application of the
FMCT to the local PO bifurcations for 
more general
Hamiltonian systems with continuous symmetries.
Section~\ref{sec4:ispm} presents
the specific application of the ISPM trace formula to the bifurcations
in the spherical RPL potential model.  Contributions of the
bifurcating orbits 
to the trace formula are discussed.  In
Sec.~\ref{sec6:comp-quant}, we compare our semiclassical results,
obtained at a bifurcation point and asymptotically far from
bifurcations, with the quantum calculations.  These results are
summarized in Sec.~\ref{sec6:concl}.  A more 
precise trace formula
based on the third-order expansions of the action integral is
presented in the Appendix. 

\section{Trace formula}
\l{sec2:trace}

The general semiclassical expression for the level density,
$g(E)=\sum_i \delta(E-E_i)$, is determined by the
energy levels $E_i$ for the single-particle Hamiltonian
$\hat{H}=\hat{T}+V$ of $\mathcal{D}$ degrees of freedom. The specific
expression generic to integrable and nonintegrable systems can be
obtained by the following trace formula in the
$2\mathcal{D}$-dimensional phase
space \cite{MFAMBpre01,MAFMptp02,MAFptp06,MFAMMSBptp99,MKApan16}:
\bea\l{pstrace}
&g_{\rm scl}(E)=\frac{1}{(2\pi\hbar)^{{\mathcal{D}}}}\Re\sum_{\rm CT}
\int \d\r^{\prime\prime} \int \d \p^{\prime}\;  
\delta \left(E - H(\r^{\prime},\p^{\prime})\right)
\nonumber\\
&\times \left|\mathcal{J}_{\rm CT}(\p^\prime_\perp,
\p^{\prime\prime}_\perp)\right|^{1/2} 
\exp\left(
\frac{i}{\hbar}\;\Phi_{\rm CT}- i 
\frac{\pi}{2} \mu^{}_{\rm CT} - i\phi^{}_{\mathcal{D}}\right).
\eea
Here $H(\r,\p)$ is the classical Hamiltonian in the 
phase-space variables $\r,\p$ and $\Phi_{\rm CT}$ is the phase integral
\bea\l{legendtrans}
\!\Phi_{\rm CT} \!&\equiv&\! 
S_{\rm CT}(\p^\prime, \p^{\prime\prime},t^{}_{\rm CT}) +
\left(\p^{\prime\prime}-\p^{\prime}\right) 
 \cdot \r^{\prime\prime} \nonumber\\
&=& S_{\rm CT}(\r^\prime, \r^{\prime\prime},E) +
\p^{\prime} \cdot \left(\r^{\prime} - \r^{\prime\prime}\right) 
\eea
(see the derivations in Ref.\ \cite{MKApan16}).  In Eq.~(\ref{pstrace}),
the sum is taken over all discrete classical trajectories (CTs) for a
particle motion from the initial point $(\r^\prime,{\bf
p}^\prime)$ to the final point $(\r^{\prime\prime},{\bf
p}^{\prime\prime})$ with a given energy $E$ \cite{MAFptp06}.  A
CT can uniquely be specified by fixing, for instance, the final
coordinate $\r^{\prime\prime}$ and the initial momentum ${\bf
p}^{\prime}$ for a given time $t^{}_{\rm CT}$ of the motion along a
CT. Here $S_{\rm CT}(\p^\prime,\p^{\prime\prime},t^{}_{\rm CT})$
is the classical action in the momentum representation,
\be\l{actionp}
S_{\rm CT}(\p^\prime,\p^{\prime\prime},t^{}_{\rm CT}) = 
-\int_{\p^\prime}^{\p^{\prime\prime}}
\d \p \cdot \r(\p)\;. 
\ee
Integration by parts relates Eq.~(\ref{actionp}) to the classical
action in the coordinate space,
\be\l{actionr}
S_{\rm CT}(\r^\prime,\r^{\prime\prime},E) = 
\int_{\r^\prime}^{\r^{\prime\prime}}
\d \r \cdot \p(\r)\;, 
\ee
by the Legendre transformation [Eq.~(\ref{legendtrans})].
The factor $\mu^{}_{\rm CT}$ is the number of
conjugate points along a CT with respect to the initial phase-space
point $(\r',\p')$. 
They are, e.g., the focal and
caustic points where the main curvatures of the energy surface 
(second derivatives of the phase integral $\Phi_{\rm CT}$)
vanish.  In addition, there are 
the turning points where these curvatures become divergent.
The number of conjugate points evaluated along a PO 
is called the Maslov index \cite{MASLOV}.
An extra phase component $\phi^{}_{\mathcal{D}}$, which is
independent of  
the individual CT,
is determined by the dimension of the 
system and the classical degeneracy $\mathcal{K}$ 
[$\phi_{\mathcal{D}}$ is zero when all orbits are isolated ($\mathcal{K}=0$), 
as defined in Ref.\ \cite{GUTZbook90}]. 

In Eq.~(\ref{pstrace}) we introduced
the local phase-space variables that 
consist of the three-dimensional 
($\mathcal{D}=3$) coordinate
$\r=\{x,y,z\}$ and momentum $\,\p=\{p_x,p_y,p_z\}$.
It is determined locally along a reference CT
so that the variables $(r_\parallel=x, p_\parallel=p_x)$
are parallel and 
$(\r_\perp=\{y,z\},\p_\perp=\{p_y,p_z\})$ are 
perpendicular to the CT \cite{GUTZbook90,SMODzp77,sclbook}. 
Here $\mathcal{J}_{\rm CT}(\p^\prime_\perp,\p^{\prime\prime}_\perp)$ 
is the Jacobian
for the transformation of the momentum component perpendicular to a
CT from the initial value $\p_\perp^\prime$ to the final value
 $\p_\perp^{\prime\prime}$. 

For calculations of the trace integral by the stationary-phase method 
(SPM), one may write
the stationary-phase conditions  
for both $\p^{\prime}$ 
and $\r^{\prime\prime}$ variables.
According to the definitions (\ref{legendtrans}) 
and (\ref{actionp}), the stationary-phase conditions are given by
\bea\l{statcond}
&&\left(\frac{\partial \Phi_{\rm CT}}{
\partial \p^{\prime}}\right)^* \equiv 
\left(\r^\prime - \r^{\prime\prime}\right)^*=0\;,\nonumber\\
&&\left(\frac{\partial \Phi_{\rm CT}}{
\partial \r^{\prime\prime}}\right)^* \equiv 
\left(\p^{\prime\prime} - \p^{\prime}\right)^*=0\;.
\eea
The asterisk indicates that quantities in  
large parentheses are
taken at the stationary point.  
Equations (\ref{statcond}) express that the
stationary-phase conditions are equivalent to the  
PO equations $(\r'',\p'')^*=(\r',{\bf
p}')^*$.  One of the SPM integrations in Eq.~(\ref{pstrace}), e.g.,
over the parallel momentum $p_{\parallel}'$ in the local Cartesian
coordinate system introduced above, is the identity because of the
energy conservation $E=H(\r'',\p'')=H(\r',{\bf
p}')$.  Therefore, it can be taken exactly.  If the system has
continuous symmetries, the integrations with respect to the
corresponding cyclic variables can be carried out exactly.  Notice
that the exact integration is performed finally also along a parallel
spatial coordinate (along the PO).

Applying the ISPM with the PO equations (\ref{statcond}), accounting
for the bifurcations and the breaking of symmetries, one may arrive at
the trace formula in terms of the sum over
POs \cite{migdalrev,sclbook,MKApan16}.  The total ISPM trace formula is
the sum over all of POs [families with the classical degeneracy ${\cal
K}\ge 1$ and isolated orbits (${\cal K}=0$)],
\be\l{dgsc}
\delta g(E) \simeq 
\delta g_{\rm scl}(E)
=\sum_{\rm PO} \delta g^{}_{\rm PO}(E)\,,
\ee
where
\begin{equation}\label{dgPO}
 \delta g^{}_{\rm PO}(E)\!=\! \Re 
\!\left\{\!A_{\rm PO} 
\exp\left[\frac{i}{\hbar} S_{\rm PO}(E) \!-\!
\frac{i\pi}{2}\mu^{}_{\rm PO} -i\phi^{}_{\mathcal{D}}\right]\right\}.
\end{equation}
The amplitude $A_{\rm PO}$ depends on the classical degeneracy
$\mathcal{K}$ and the stability of the PO. In the exponent phase,
 $S_{\rm PO}(E)=\oint\p\cdot d\r$ is the action and
$\mu^{}_{\rm PO}$ is the Maslov index
\cite{GUTZbook90,SMepp76,sclbook,migdalrev,MKApan16}.

\section{Fedoriuk-Maslov catastrophe theory} 
\l{sec3:cat3}

In this section we present the essence of the ISPM, following
basically the Fedoriuk-Maslov catastrophe theory 
(see Refs.\ \cite{FEDoryuk_pr,MASLOV,FEDoryuk_book1}). 

\subsection{Caustic and turning points}
\l{sec31:caustics}

Let us assume that the integration interval in one of the integrals of
Eq.~(\ref{pstrace}) over a phase-space variable $\x$ contains a
stationary catastrophe point where the second derivative of the phase
integral $\Phi$ is zero [see Eq.~(\ref{statcond})]. This catastrophe
integral $\mathcal{I}(\kappa,\epsilon)$ can be considered as a
function of the two dimensionless parameters
\begin{equation} 
\mathcal{I}(\kappa,\epsilon) = \int_{\x_-}^{\x_+} \d \x 
\;A(\x,\epsilon)\,\exp\left[i\,\kappa\, \Phi(\x,\epsilon)\right], 
\label{catint} 
\end{equation} 
where $A(\x,\epsilon)$ 
is  
the amplitude and $\Phi(\x,\epsilon)$ the 
dimensionless phase 
integral, which is proportional to $\Phi_{\rm CT}$ given by 
Eq.~(\ref{pstrace}).
One of these parameters $\kappa$ is related to a large
semiclassical parameter
$\kappa \propto 1/\hbar \rightarrow \infty $,
when $\hbar \rightarrow 0$, through the relationship
$\kappa \Phi=\Phi_{\rm CT}/\hbar$
(see also the Appendix 
for a  
clear example).
Another critical parameter $\epsilon$ is a small dimensionless  perturbation 
of the phase integral $\Phi(\x,\epsilon)$ 
[Eq.~(\ref{legendtrans})] and the amplitude $A(\x,\epsilon)$
through the potential $V(\r,\epsilon)$.  
For instance, $\epsilon$ can be a dimensionless distance from
the catastrophe 
point by perturbing 
the parameter of a potential  $V(\r,\epsilon)$
\cite{migdalrev}, e.g., 
the  deformation and diffuseness parameters
(see examples in Refs.~\cite{sclbook,migdalrev,MKApan16}).  
In Eq.~(\ref{catint}), the integration limits $\x_\pm$
crossing the catastrophe
point $\x^\ast(0)$ are
generally assumed to be finite.

We assume also that the integral (\ref{catint}) has 
the simplest (first-order) caustic-catastrophe point 
$\x^*(\epsilon)$ at $\epsilon=0$ defined by\footnote{
In general, the caustic point of the $n$th order $(n\geq 1)$
is defined as the point where the
derivatives up to the ($n+1$)th order vanish but the $(n+2)$th
derivatives remain finite.}
\bea\l{caustpoint} 
&\Phi'\left(\xi^\ast\right) = 
0,\qquad  \Phi''\left(\xi^\ast\right) = 0,\nonumber\\ 
&\Phi'''\left(\xi^\ast\right) = 
 O\left(\epsilon^0\right) \quad \mbox{at}\quad \x^\ast=\x_0=\x^\ast(0),
\eea
where the asterisk indicates
that the derivatives with respect to $\xi$ are taken at $\x=\x^\ast$. 
The mixed derivative $(\partial^2 \Phi/\partial \x\partial \epsilon)^\ast$,
$A^\ast=A(\xi^\ast,\epsilon)$,
is assumed to be of the zeroth order in $\epsilon$ 
as well as the third derivative
in Eq.~(\ref{caustpoint}).
In the limit $\epsilon \rightarrow 0$ at large 
$\kappa$ for the caustics,
the two simple stationary points $\x^\ast(\epsilon)$ coincide and form 
one caustic point given by Eq.~(\ref{caustpoint}). 
To remove the indetermination,
let us consider a small perturbation of the catastrophe integral 
$I(\epsilon,\kappa)$ [Eq.~(\ref{catint})] 
by changing $\epsilon$ through the phase 
function $\Phi(\x,\epsilon)$ and the amplitude $A(\x,\epsilon)$ 
near the caustic point $\epsilon=0$. 
For any small nonzero $\epsilon$ we first study the 
expansion of the function
$\Phi(\x,\epsilon)$ over $\epsilon$ in 
a power series near the stationary point $\x^*(\epsilon)$ for 
a small $\epsilon$,
\bea\l{Sexp}
\Phi(\x,\epsilon) &=& \Phi^\ast 
+ \frac{1}{2} \,\Phi''\left(\xi^\ast\right) 
\,(\x-\x^\ast)^2 \nonumber\\
&+&\frac{1}{6}\,\Phi'''\left(\xi^\ast\right)\,(\x-\x^\ast)^3 + 
\cdots~. 
\eea 
Similarly, for the amplitude
expansion, one has
\be\l{Aexp}
A(\x,\epsilon) = A^\ast + A'\left(\x^\ast\right)\;
(\x-\x^\ast) + \cdots~.
\ee
The asterisks in  
Eqs.~(\ref{Sexp}) and (\ref{Aexp}) 
indicate that 
the derivatives with respect to $\x$ are taken at $\x=\x^*(\epsilon)$
for a small but finite $\epsilon$, 
$\Phi^* =\Phi(\x^*,\epsilon)$. 
Using a small perturbation of the action $\Phi(\x,\epsilon)$ and amplitude 
$A(\x,\epsilon)$ by $\epsilon$ variations, 
one finds the first derivative in Eq.~(\ref{caustpoint}) and
the second derivative in Eq.~(\ref{Sexp}) as small but 
nonzero quantities. 
For asymptotic values of $\kappa \rightarrow \infty$, one may 
truncate
the series (\ref{Sexp}) for the phase $\Phi(\x,\epsilon)$ 
in the exponent of the catastrophe integral $\mathcal{I}(\kappa,\epsilon)$
[Eq.~(\ref{catint})]
and corresponding one (\ref{Aexp}) for its amplitude $A(\x,\epsilon)$ 
at the third and 
zeroth orders, respectively, 
keeping, however, a small nonzero $\epsilon$.

In the following, defining $\xi$ as a dimensionless variable,
one can simply take the second derivative of a phase $\Phi$ in
Eq.~(\ref{Sexp}), divided by 2, as the parameter $\epsilon$, 
\be\l{eps}
\epsilon=\frac12\Phi''\left(\x^\ast\right)=\sigma|\epsilon|\;,
\quad \sigma=\sgn(\epsilon)=\sgn[\Phi''(\x^*)]\;.
\ee
By definition of the simplest caustic point of the first order 
[Eq.~(\ref{caustpoint})],
the third derivative of the phase $\Phi$  
near the caustic point $\epsilon=0$ is not zero
at any small $\epsilon$.
Therefore, one can truncate the expansion of the phase
integral (\ref{Sexp}) up to the third order as
\be\l{norm3}
\Phi=\Phi^\ast + \epsilon\;(\x-\x^*)^2
     + a\;(\x-\x^*)^3\;,
\ee
where 
\be\l{aQ}
a=\frac16 \Phi'''\left(\xi^\ast\right)\;.
\ee
Here we assume that $A(\x^*,\epsilon)$
has a finite nonzero
limit at $\epsilon \to 0$.  Therefore, it can be cut at zeroth
order in $\x-\x^*$, namely, 
 $A(\xi,\epsilon) \approx  A^*=A(\xi^*,0)$.

Following Refs.~\cite{FEDoryuk_pr,FEDoryuk_book1},
one can consider the linear transformation 
of the coordinate from $\x$ to $z$ as
\begin{equation} 
\xi-\xi^\ast =  \Upsilon + \Lambda\,z\;, \quad
\label{transform1} 
\end{equation}
where 
\be\l{transformcoeff}
\Upsilon=-\frac{\epsilon}{3a}\;,\quad \Lambda=
\frac{1}{\left(3a\;\kappa\right)^{1/3}}\;.
\ee
Notice that both 
coefficients of this linear transformation,
$\Upsilon(\epsilon)$ and $\Lambda(\kappa)$, depend on different
critical parameters $\epsilon$ and $\kappa$, respectively
[$\Upsilon(0)=\Lambda(\infty)=0$].  Substituting
Eq.~(\ref{transform1}) into Eq.~(\ref{norm3}), one can express the
catastrophe integral (\ref{catint}) in an analytical form,
\bea\l{catint2ai} 
&\mathcal{I}(\kappa,\epsilon)= \pi \Lambda A^\ast\, 
\exp\left(i \,\kappa \,\Phi^* + \frac{2i\sigma}{3}\, 
w^{3/2}\right)\,\nonumber\\ 
&\times\left[\mbox{Ai}(-w,\mathcal{Z}_{-},\mathcal{Z}_{+}) + 
i\,\mbox{Gi}(-w,\mathcal{Z}_{-},\mathcal{Z}_{+})\right]\;. 
\eea 
The generalized incomplete Airy and Gairy integrals are defined in a
similar way as the standard ones but with the finite integration
limits,
\be\l{airynoncomp} 
\begin{Bmatrix} \mbox{Ai} \\ \mbox{Gi}\end{Bmatrix}(-w,z_1,z_2) 
\!=\!\frac{1}{\pi}\int_{z_{1}}^{z_{2}}\;\d z
\begin{Bmatrix} \cos \\  \sin \end{Bmatrix}
\left(-w\;z\!+\!\frac{z^3}{3}\right).
\ee
The argument $w$ of these functions and finite integration 
limits $\mathcal{Z}_\pm$ in 
Eqs.~(\ref{catint2ai}) and (\ref{airynoncomp}) are given by 
\be\l{w} 
w= \frac{\kappa^{2/3}\;\epsilon^2}{\left(3a\right)^{4/3}} > 0,
\ee
and
\be\l{zpm}
\mathcal{Z}_\pm = \frac{\xi_\pm -\xi^*}{\Lambda} +\sigma
\sqrt{w}\;.
\ee
As can be seen from a cubic form of the phase in parenthesis on 
the right-hand side of Eq.~(\ref{airynoncomp}),
the caustic catastrophe can be considered as a crossing point of  
the two simple close stationary-point curves for any 
small nonzero $\epsilon$,
\be\l{statpoints3}
z_\pm^*(\epsilon)=\pm \sqrt{w(\epsilon)}\;.
\ee 
They degenerate into one caustic point
$z_\pm^*(\epsilon) \rightarrow 0$ (\ref{caustpoint}) in the limit 
$\epsilon \rightarrow 0$ because of $w \rightarrow 0$ at any finite $\kappa$. 
The final result is a sum of the contributions of these 
stationary points. 
Note that, according to Eq.~(\ref{w}),
the value of $w$ is large for any nonzero $\epsilon$ when 
$\kappa $ becomes large.  Nevertheless, it becomes 
small for a large fixed 
finite nonzero $\kappa$ when $\epsilon$ is small. 
However, we 
may consider 
both cases by using the same formula 
(\ref{catint2ai}) because
the two parameters $\epsilon$ and $\kappa$ appear in 
(\ref{w}) through one parameter $w$ for a finite constant $a$. 

\subsection{The Maslov theorem}
\l{sec32:maslov}

\begin{figure}
\centering
\includegraphics[scale=.4]{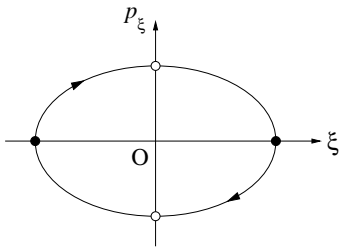}
\caption{\label{fig1}
Illustration of the phase-space flow around a PO at the origin
in the case of a regular trajectory on the torus.
Open and  closed circles  
represent caustic points
($\partial p_\xi/\partial\xi=0$) and turning points
($\partial p_\xi/\partial\xi=\pm\infty$), respectively.
The trajectory around the PO turns clockwise in the $(\xi,p_\xi)$ plane
and the curvature $\Phi''\propto \partial p_\xi/\partial\xi$ changes its
sign from positive to negative when the trajectory crosses the caustic point.
At the crossings of the turning points, the curvature changes from
$-\infty$ to $+\infty$.}
\end{figure}

For any small nonzero $\epsilon$, one may find a value of
$\kappa$ for which $w$ is so large that
the two above-mentioned stationary points $z_\pm^*(\epsilon)$ can be  
treated separately in the SPM.
At the same time, $w$ can be sufficiently 
small 
that the stationary
points $|z_\pm^*|$ are much smaller than the integration limits
$|\mathcal{Z}_\pm|$.
In practice, it is enough to consider $\mathcal{Z}_{+}>0$
and $\mathcal{Z}_{-}<0$, and 
we split the integration interval into two parts, i.e., from
$\mathcal{Z}_{-}$ to $0$ and from $0$ to $\mathcal{Z}_{+}$,
in order to separate the contributions of
negative and positive stationary points
$-\sqrt{w}$ and $\sqrt{w}$.
According to the phase-space flow
around the PO (see Fig.~\ref{fig1}),
the curvature $\Phi''\propto \partial p_{\x}/\partial \x$ ($p_{\x}$
denotes the momentum conjugate to $\x$) will always change its sign
from positive to negative
at the caustic catastrophe point \cite{Ozoriobook}.
Let us consider such crossings with a catastrophe point.
Outside  
the catastrophe point, where
$\epsilon>0\; (\sigma=1)$, 
the stationary point $\x^*$ corresponds to
$z_+^*$ and its contribution is given by the
second integral over positive $z$.
For this integral, one can extend the upper integration limit
$\mathcal{Z}_{+}$ to 
infinity, like in the SSPM, because  
$\mathcal{Z}_{+} \propto \kappa^{1/3} \gg 1$ and 
$\mathcal{Z}_{+} \gg z_{+}^*=\sqrt{w} \propto \kappa^{1/3} \epsilon$ for 
a small finite $\epsilon$ [see Eq.~(\ref{w})].
Within this approximation, one can use
the standard complete Airy and Gairy functions
\begin{equation} 
\begin{Bmatrix}\mbox{Ai} \\ \mbox{Gi}\end{Bmatrix}(-w)
=\frac{1}{\pi}\int_{0}^{\infty}\d z
 \begin{Bmatrix} \cos \\ \sin \end{Bmatrix}
 \left(-w\,z+\frac{z^3}{3}\right). \label{airystandard}
\end{equation}
They correspond to the limits 
$z_1=\mathcal{Z}_{-}=0$ and $z_2=\mathcal{Z}_{+}=\infty$
in Eq.~(\ref{airynoncomp}). 
Using the asymptotic form of the Airy and Gairy functions for
$w \rightarrow \infty$,
\be\l{AiGiAsympt}
\begin{Bmatrix}\mbox{Ai} \\ \mbox{Gi}\end{Bmatrix}(-w)
 \to \frac{1}{\sqrt{\pi}w^{1/4}}
 \begin{Bmatrix}\sin \\ \cos\end{Bmatrix}
 \left(\frac23w^{3/2}+\frac{\pi}{4}\right)\;,
\ee
one can evaluate the contribution $\mathcal{I}_+$ 
of the positive stationary point to Eq.~(\ref{catint2ai}),
asymptotically 
far from the caustic point (\ref{caustpoint}).
Then one obtains the same result 
as one 
would get by the standard second-order expansion of the phase $\Phi$ 
(and zeroth order of the amplitude $A$)
at a simple stationary point $\x^*$,
\begin{gather}
\mathcal{I}_+(\kappa,\epsilon)
\simeq \pi\Lambda A^*\exp\left(i\kappa\Phi^*+ \frac{2i}{3} w^{3/2}\right)
\left[\mbox{Ai}(-w)\right. \nonumber \\
+\left. i\mbox{Gi}(-w)\right]
\rightarrow 
\sqrt{\frac{2\pi}{\kappa|\Phi''(\x^*)|}}A^*\exp\left(
{i\kappa\Phi^*+\frac{i\pi}{4}}\right)\;.
\label{eq:int_before}
\end{gather}
On the other side of the
crossing with the catastrophe point, $\epsilon<0\; (\sigma=-1)$,
the stationary point $\x^*$ corresponds
to $z_-^*$.  Therefore, one should consider the 
other part of the integral
$\mathcal{I}_-$ over the negative values of $z$.
Obviously, considering
in an analogous way with a change of the integration variable $z \to -z$,
one obtains
\begin{gather}
\mathcal{I}_-(\kappa,\epsilon)
\simeq \pi\Lambda A^* \exp\left(i\kappa\Phi^*-\frac{2i}{3}w^{3/2}\right)
\left[\mbox{Ai}(-w)\right. \nonumber \\
-\left.i\mbox{Gi}(-w)\right]
\rightarrow \sqrt{\frac{2\pi}{\kappa|\Phi''(\x^*)|}}A^*
\exp\left(i\kappa\Phi^*-\frac{i\pi}{4}\right)\;.
\label{eq:int_after}
\end{gather}
Comparing the rightmost expression in
Eq.~(\ref{eq:int_after})
with that in Eq.~(\ref{eq:int_before}), one sees a
shift of phase by $-\pi/2$.
Thus, the famous Maslov theorem \cite{MASLOV}
on the shift of the phase $\Phi$
by $-\pi/2$ at each simple caustic point 
(\ref{caustpoint}) of the CT 
(in particular, the PO) in the SSPM 
(i.e., the Maslov index is increased by one)
has been proved by using the Fedoriuk 
catastrophe method \cite{FEDoryuk_pr}.

For the case of a turning point, one 
has the conditions (\ref{caustpoint}) with only the 
replacement of zero by infinity
in the second derivative. In this case, Fedoriuk used a linear 
coordinate transformation from $z$ to $\widetilde{z}$,
which has the form
\be\l{turningtrans}
z=\frac{\widetilde{\Lambda}}{\epsilon}\widetilde{z} + \widetilde{\Upsilon},
\ee
where $\widetilde{\Lambda}$
and $\widetilde{\Upsilon}$ are new constants 
that are not
singular in $\epsilon$ ($\widetilde{\Lambda}$ is independent of 
$\epsilon$).
This transformation reduces the 
turning-point singularity to the 
caustic-point one. Indeed,
the divergent second derivative of the phase $\Phi$ over the variable $z$
is transformed, in a new variable, to its zero value
\be\l{newderPhi}
\frac12\left(\frac{\partial^2\Phi}{\partial^2\widetilde{z}}\right)^*=
\frac{\widetilde{\Lambda}^2}{\epsilon}\;
\ee
(see the second paper in Ref.~\cite{FEDoryuk_pr}). 
Therefore, one obtains the same shift by $-\pi/2$ at each simplest 
turning point
(per a sign change of one momentum component perpendicular to
the boundary) along a CT.  Finally, 
the next part of the Maslov theorem \cite{MASLOV}
concerning the Maslov index generated by such a 
turning  point has been proved, too, within the same catastrophe 
theory of Fedoriuk \cite{FEDoryuk_pr}.

\section{Symmetry-breaking and bifurcations}
\l{sec4:bifurcations}

Assuming a convergence of the expansion
[Eq.~(\ref{Sexp})] to the second order, as shown in the 
preceding section, one can use the same approach 
within the simplest ISPM2 of second order\footnote{We call the 
$n$th-order ISPM  
the ISPM$n$, in which we use
the expansion of the phase 
integral $\Phi$ up to the $n$th-order terms and the amplitude up to 
the $(n-2)$th order near the stationary point. 
The simplest ISPM for $n=2$ is called  
the ISPM2.}
to investigate the symmetry-breaking and bifurcation problems in the
POT.  Therefore, one arrives at a sum over the separate contributions
of different kinds of 
isolated and degenerated orbits, as in the
derivation of the Maslov theorem but within the finite
integration limits. The latter is important because the bifurcation
point is located at a boundary of the classically accessible region.
This section is devoted to the extension of the FMCT to the
bifurcation catastrophe problems.

In the presence
of continuous symmetries, the stationary points form a family of POs 
that cover a 
$(\mathcal{K}$$+$$1)$-dimensional submanifold 
$\mathcal{Q}_{\rm PO}$ of phase space, whereby
$\mathcal{K}$ is the classical degeneracy of the PO family.  
The integration over $\mathcal{Q}_{\rm PO}$ must be performed exactly.
In any 
systems with continuous symmetries,
it is an advantage to transform the phase-space 
variables 
from the Cartesian to the corresponding action-angle variables (see, e.g., 
Ref.~\cite{BT76,CREAGHpra97}).
Then the action $\Phi_{\rm CT}$ in Eq.~(\ref{pstrace})
is independent of the angle variables conjugate to
the conserving action variables
and the integrations over these cyclic angle variables are 
exactly carried out.
For the integrable case, for instance, integrating over the 
remaining action variables and using the 
standard SPM, one obtains the so-called Berry-Tabor trace formula \cite{BT76}.
Under the existence of additional symmetries  
such as SU(3) or O(4), 
some of the integrations over the action variables
can also be  
performed exactly because of a higher degeneracy.
For partially integrable systems, the integrations over partial set of 
cyclic variables also  
greatly simplify the ISPM derivations of the trace 
formula (\ref{pstrace}) near the bifurcations.

To solve the bifurcation problems,
some of the SPM integrations have to be done
in a more exact way.
For definiteness, we  
will consider first a simple bifurcation 
defined
as a caustic point of the first order [Eqs.~(\ref{Sexp}) and (\ref{Aexp})]
where the degeneracy parameter $\mathcal{K}$ is locally 
           increased by one.
In the SPM, after performing  
exact integrations over 
a submanifold $\mathcal{Q}_{\rm PO}$, one uses an expansion of
the action phase $\Phi_{\rm CT}$ in phase space variables
$\xi=\{\r'',\p'\}_\perp$, perpendicular to $\mathcal{Q}_{\rm PO}$ in the
integrand of Eq.~(\ref{pstrace}) over $\xi$ near the stationary point
$\xi^\ast$, 
\bea\l{ispmexp}
\!\!&\Phi_{\rm CT}(\xi)=\Phi^{}_{\rm PO} + \frac12\Phi_{\rm PO}''(\xi^\ast) 
(\xi-\xi^{\ast})^2 \nonumber\\
&+ \frac16\Phi_{\rm PO}'''(\xi^\ast) (\xi-\xi^{\ast})^3+
\cdot\cdot\cdot\;,
\eea
where
\be\l{PhiPO}
\Phi_{\rm PO}=\Phi^\ast_{\rm CT}=\Phi_{\rm CT}(\xi^\ast)\;,
\ee
$\xi^\ast$ is the stationary point, $\xi^\ast=\xi^{}_{\rm PO}$, and
$\Phi^\ast_{\rm CT}=\Phi_{\rm CT}\left(\xi^\ast\right)=\Phi_{\rm PO}$.
To demonstrate the key point of our derivations of the trace formula, 
we focus on one of the phase-space variables in Eq.~(\ref{pstrace}),
denoted by $\xi$, which is associated with a
catastrophe behavior.
In the standard SPM, the above expansion is truncated at the second-order term
and the integration over the variable $\xi$ is extended to $\pm\infty$.
The integration can be performed analytically and yields a Fresnel integral
(see, e.g., Refs.~\cite{sclbook,migdalrev}). 

However, one encounters a singularity in the SSPM that 
is related to the zero or  
infinite value of $\Phi_{\rm PO}''(\xi^\ast)$, 
while $\Phi_{\rm PO}'''(\xi^\ast)$ remains finite in the simplest case
under consideration. 
This singularity occurs when a PO (isolated or degenerated)
undergoes a simple bifurcation at 
the stationary point $\xi^\ast$ under the 
variation of a parameter of the potential 
(e.g., energy, deformation, or 
surface diffuseness).
The SSPM approximation to the Fresnel (error) functions
by the Gaussian integrals breaks down
because one has a divergence.

Notice that the bifurcation problem is similar to the 
caustic singularity considered by Fedoriuk within the catastrophe theory 
(see Sec.~\ref{sec3:cat3} and Refs.~\cite{FEDoryuk_pr,FEDoryuk_book1}).  
The FMCT is adopted, however, 
for the specific position of such a singularity at the end point in
the phase-space volume accessible for a classical motion 
(see the Introduction and also Ref.\ \cite{MAFptp06}).

In systems with continuous symmetries, the 
orbit at the end point causes the bifurcation where it coincides with
one of the rational tori that 
appears in the transition from the unphysical to the physical region. 
The contribution of this end-point orbit is derived
using a local phase-space variable $\xi$ along it, independently of the
torus orbits. 
Near the bifurcation point, the contribution of these
end-point orbits is mostly included 
to the newborn orbit term.
Therefore, we should consider a kind of separation of the phase
space occupied by the newborn and 
end-point orbits to evaluate their contributions to the trace integral
near the bifurcation point. Below we will call 
the latter an end-point manifold. 
By definition of the end-point manifold,
its measure is zero at the 
bifurcation limit where the minimal $\x_{-}$ and maximal $\x_{+}$
coincide with the stationary point $\x^\ast$, i.e.,
\be\l{bifboundlim}
\x_{-} \rightarrow \x_{+} \rightarrow \x^\ast
\quad
(\Phi_{\rm PO}''\to 0).
\ee
Thus, although the contributions of POs participating in
the bifurcation
are considered separately, the parent orbit contribution vanishes at
the bifurcation point and 
there no risk of double counting.
To describe the transition from the bifurcation point to the
asymptotic region, one should properly define the end-point
manifold.
We need it to extract the
additional contribution by the end-point orbit that is
not covered by the term for the newborn orbit.
The detailed treatment of such a transition 
is still open, but is beyond the scope of the present paper.

We are ready now to employ what we call the improved 
stationary-phase method  
\cite{MFAMMSBptp99,MFAMBpre01,MAFMptp02,MAFptp06}
evaluating the trace integral for the semiclassical level density.
Hereby, the integration over $\xi$ in Eq.~(\ref{pstrace}) is restricted 
to the finite limits defined by the classically allowed phase-space 
region through the energy-conserving 
$\delta$ function  
in the integrand of Eq.~(\ref{pstrace}).
The phase and 
 amplitude are expanded around the stationary point
up to the second- and zeroth-order terms in $\xi -\xi^*$,
respectively, and to higher-order terms if necessary.

In the simplest version of the ISPM
(ISPM2), the expansion of the phase is truncated
at second order, keeping the finite integration limits $\xi_{-}$ and $\xi_{+}$ 
given by the accessible region of the classical motion in 
Eq.~(\ref{pstrace}).
It will lead to a factor like\footnote{
For the case of several variables $\xi$ for which we find find zeros and 
infinities in eigenvalues of the matrix
with second-order derivatives of $\Phi_{\rm PO}(\xi)$ at $\x=\xi^*$, we 
diagonalize this matrix and reduce the Fresnel-like integrals to products of 
error functions similar to Eq.~(\ref{errfuns}).}
\bea\l{errfuns}
&e^{i\Phi_{\rm PO}/\hbar}
\int_{\xi_{-}}^{\xi_{+}} \exp\left[\frac{i}{2\hbar}\Phi_{\rm PO}''\; 
(\xi-\xi^{\ast})^2\right]
\mbox{d} \xi \nonumber\\
&\propto \frac{1}{
\sqrt{\Phi_{\rm PO}''}}\; 
e^{i\Phi_{\rm PO}/\hbar}\;
\erf\left[\mathcal{Z}_{-},\mathcal{Z}_{+}\right],
\eea
where $\erf(z_1,z_2)$ is the generalized error function,
\be\l{ispmlimits}
\erf(z_1,z_2)=\frac{2}{\sqrt{\pi}}\int_{z_1}^{z_2}e^{-z^2}{\mbox d}z
=\erf(z_2)-\erf(z_1), 
\ee
with the complex arguments
\be\l{comparg}
\mathcal{Z}_{\pm}=\left(\xi_{\pm}-\xi^*\right)
\sqrt{-\frac{i}{2\hbar}\Phi_{\rm PO}''}\;.
\ee
Note that the
expression (\ref{errfuns}) has no divergence at the 
bifurcation point where $\Phi_{\rm PO}''(\xi^\ast)=0$, since the error
function (\ref{ispmlimits}) also goes to zero as
\be
\erf(\mathcal{Z}_-,\mathcal{Z}_+)
\propto \mathcal{Z}_+-\mathcal{Z}_-
\propto(\x_+-\x_-)\sqrt{\Phi''_{\rm PO}}\;.
\ee
Thus, the factor $\sqrt{\Phi_{\rm PO}''}$ in the denominator
of Eq.~(\ref{errfuns})
is canceled with the same in the numerator.
In addition, also taking into account
Eq.~(\ref{bifboundlim}) and the discussion around it, one finds the
zero contribution of
the end-point term in
the trace formula at the bifurcation point,
which seems to be a consistent semiclassical picture.

This procedure is proved to be valid in the semiclassical limit 
$\kappa \rightarrow \infty$
by the FMCT \cite{FEDoryuk_pr,MASLOV,FEDoryuk_book1}.
In this way, we can derive the separate PO contributions that 
are free of divergences, discontinuities, and double counting
at any bifurcation point.
The oscillating part of the
level density can be approximated by the semiclassical 
trace formula (\ref{dgsc}).
In Eq.~(\ref{dgsc}), the sum runs over all periodic orbits 
(isolated or degenerated) in
the classical system. The term
$S_{\rm PO}(E)$ is the action integral along a PO. 
The amplitude $A_{\rm PO}(E)$ (which, in general, is complex) is 
of the order of the phase-space volume occupied by CTs. 
The factor given in Eq.~(\ref{errfuns}) depends on the degeneracies
and stabilities of the POs, respectively 
(see Sec.~\ref{sec3:cat3}).

Notice that any additional exact integration
in Eq.~(\ref{pstrace}) with respect to
a bifurcation (catastrophe) variable 
of the improved SPM can lead to an enhancement of the amplitude $A_{\rm PO}$
in the transition
from the bifurcation point to the asymptotic region.
This enhancement is of the order $1/\hbar^{1/2}$ as compared to the 
result of the standard SPM integration. 
In particular, for the newborn family with
the extra degeneracy $\Delta\mathcal{K}$ higher than 
that of the parent PO, one has such enhancements
of the order $1/\hbar^{\Delta\mathcal{K}/2}$ near the bifurcation.

The trace formula (\ref{dgsc}) thus relates the quantum oscillations 
in the level 
density to quantities  
that are purely determined by the classical system.
Therefore, one can understand the shell effects
in terms of classical pictures. 
The sum over POs in Eq.~(\ref{dgsc}) is asymptotically 
correct to the leading order in $1/\hbar^{1/2}$ and 
it is hampered by convergence problems \cite{GUTZbook90}.
However, one is free from those problems by taking
the coarse-grained level density
\be\l{deltadenav}
\delta g^{}_{\Gamma,{\rm scl}}(E)=
\sum_{\rm PO} \delta g^{\rm scl}_{\rm PO}(E)\,
\exp\left\{-\left(\frac{\Gamma t^{}_{\rm PO}}{2\hbar}\right)^{\!2}\right\},
\ee
where $\Gamma$ is an averaging width. The value of $\Gamma$ 
is sufficiently smaller than
the distance between the major shells 
near the Fermi surface (see 
Refs.~\cite{BBann72,SMepp76,SMODzp77,sclbook,migdalrev,MKApan16}). Here
$t^{}_{\rm PO}=\partial S_{\rm PO}(E)/\partial E$
is the period of the
particle motion along a PO
taking into account its repetition number.
We see that, depending on the 
smoothing width $\Gamma$,
longer orbits are automatically 
suppressed in the above expressions and the PO sum converges, which it 
usually does not \cite{GUTZbook90} for nonintegrable systems in the limit 
$\Gamma \rightarrow 0$.
Thus, one can
highlight the major-shell structure in the level density
using a smoothing width  
that is much larger than the mean 
single-particle level spacing but 
smaller than the main shell spacing (the distance between major shells) near 
the Fermi surface. Alternatively, a finer shell structure can
be considered by using essentially
smaller smoothing widths,  
which is important for
studying the symmetry-breaking (bifurcation) phenomenon
associated with longer POs. 

It is an advantage of this approach that the
major-shell effects in $\delta g(E)$ can 
often be explained semiclassically 
in terms of only a few of the shortest POs in the system. 
Examples will be given in Sec.~\ref{sec4:ispm}.
However, if one wants to study a finer shell
structure, specifically
at large deformations, some longer orbits have to be included.
Hereby, bifurcations of POs
play a crucial role, as it will be exemplified in Sec.~\ref{sec4:ispm}.

\section{ISPM for the spherical RPL potential model}
\l{sec4:ispm}

Let us apply the general FMCT to the RPL potential
as an analytically solvable example.

\subsection{Scaling property}
\l{sec5A:scal}

A realistic mean-field potential for nuclei and metallic clusters
is given by the well-known Woods-Saxon (WS) potential
\be
V_{\rm WS}(r)=-\frac{W_0}{1+\exp\frac{r-R}{a}}\;,
\ee
where $W_0$ is the depth of the potential,
$R$ is the nuclear radius, and $a$ is the surface diffuseness.
As suggested in Refs.~\cite{KAijmpe04,KAprc12}, this potential
can be approximated by the RPL potential
for a wide range of mass numbers as
\be\l{ramod}
V_{\rm WS}(r) \approx -W_0 + \frac{W_0}{2}(r/R)^\alpha\;,
\ee
with an appropriate choice of the radial power parameter $\alpha$.
One finds  
good agreement of the quantum
spectra for the approximation (\ref{ramod}) to the WS potential
up to and around the Fermi energy $E^{}_F$.
Eliminating the constant term  
on the right-hand side of
(\ref{ramod}), we define the RPL model Hamiltonian as
\be\l{potenra}
H=\frac{p^2}{2m}+V_0(r/R_0)^\alpha, \quad V_0=\frac{\hbar^2}{mR_0^2}\;,
\ee
where $m$ is the mass of a particle and $R_0$ is an arbitrary length
parameter.

In the spherical RPL model (\ref{potenra}),
as well as in general spherical potential models,
one has the diameter and circle POs  
that form the
two-parameter ($\mathcal{K}=2$) families.
The diameter and circle POs have minimum and maximum values of the
angular momentum $L=0$ and $L_C$, respectively.
They correspond to the end points
of the energy surface $H(I_r,L)=E$ implicitly given by the
relationship (see, e.g., Ref.~\cite{MKApan16})
\be\l{ensurf}  
I_r=\frac{1}{\pi}\int_{r_{\rm min}}^{r_{\rm max}} p_r \d r= I_r(E,L), 
\ee
where $p_r$ is the radial momentum
\be\l{pr}
p_r=\sqrt{p^2(E,L)-L^2/r^2},
\ee 
with $p$ the particle momentum
\be\l{p}
p(E,L)=\sqrt{2m\left[E-V(r)\right]}\;.
\ee
The integration limits $r_{\rm min}$ and $r_{\rm max}$
in Eq.~(\ref{ensurf}) are
functions of the energy $E$ and angular momentum $L$
for a given spherical potential $V(r)$.  In the spherical RPL potential,
they are defined by the two real roots of the transcendent equation
for the variable $r$:
\be\l{pr2eq}
 p_r(L,\alpha)
 \equiv \sqrt{2m\left(E-V_0(r/R_0)^\alpha\right)-L^2/r^2}=0\;.
\ee
Another key quantity in the POT
is the curvature of the energy surface (\ref{ensurf})
\be\l{curvature}
K=\frac{\partial^2 I_r}{\partial L^2}\;.
\ee

Using the invariance of the equations of motion
under the scale transformation
\bea\l{scaling}
&&\r\to s^{1/\alpha}\r, \quad
\p\to s^{1/2}\p\;, \nonumber\\
&&t \to s^{1/2-1/\alpha}t \quad
\mbox{for} \quad E\to sE\;, 
\eea
one may factorize 
the action integral
$S_{\rm PO}(E)$ along the PO as
\bea\l{actionsc}
&&S_{\rm PO}(E)=\oint_{{\rm PO}(E)} \p \cdot \d \r \nonumber \\
&&=\left(\tfrac{E}{V_0}\right)^{\frac12+\frac{1}{\alpha}}
\oint_{{\rm PO}(E=V_0)} \p \cdot \d \r
\equiv\hbar\epsi\tau^{}_{\rm PO}\;. 
\eea
In Eq.~(\ref{actionsc}) we define the dimensionless
variables $\epsi$ and $\tau^{}_{\rm PO}$ 
as classical characteristics of the particle motion
\be\l{eq:scaledentau1}
\epsi=\left(E/V_0\right)^{\frac12+1/\alpha}
\ee
and
\be\l{eq:scaledentau2}
\tau^{}_{\rm PO}=\frac{1}{\hbar}\oint_{{\rm PO}(E=V_0)}
\p \cdot \d \r\;.
\ee
We call them the
scaled energy and the scaled period, respectively.
To realize the advantage of the scaling
invariance under
the transformation (\ref{scaling}),
it is helpful to use
$\epsi$ and $\tau^{}_{\rm PO}$ in place of the
energy $E$ and the period $t_{\rm PO}$ for the particle motion along
a PO, respectively.
In the harmonic-oscillator (HO) limit ($\alpha\to 2$), $\epsi$
and $\tau^{}_{\rm PO}$ are 
proportional to 
$E$ and $t^{}_{\rm PO}$,
while in the cavity limit
($\alpha \rightarrow \infty$), they are proportional 
to the momentum $p$ and the geometrical PO 
length $\mathcal{L}_{\rm PO}$,
respectively.

\begin{figure}
\begin{center}
\includegraphics[width=\linewidth,clip=true]{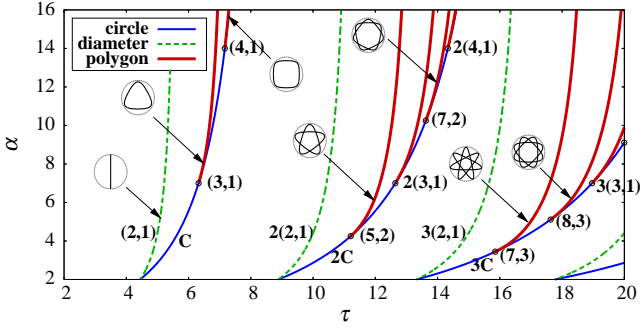}
\end{center}
\caption{Scaled periods $\tau^{}_{\rm PO}$ (horizontal axis) of short
periodic orbits PO plotted as functions of the power parameter
$\alpha$ (vertical axis).  Thin solid blue curves are circle orbits $MC$,
dashed green curves are diameter orbits $M$(2,1), and thick solid red
curves are polygonlike orbits $M(n_r,n_\theta)$ $(n_r>2n_\theta)$ 
which bifurcate from the circle orbits $MC$ at the bifurcation points
indicated by open circles.  
}
\label{fig2}
\end{figure}

The PO condition (\ref{statcond})
determines 
several PO families in the RPL potential, namely, the 
polygonlike ($\mathcal{K}=3$), the circular, and the diametric
($\mathcal{K}=2$)
POs.  This condition for the integrable spherical Hamiltonian is 
identical to a resonance condition, which is expressed in
the spherical variables $r, \theta, \varphi$ as
\be\l{rescond}
\omega^{}_r/\omega^{}_\theta =n^{}_r/n^{}_\theta~,\qquad 
\omega^{}_\theta \equiv \omega^{}_\varphi\;,
\ee
where $\omega_r$, $\omega_\theta$, and $\omega_\varphi$ are
frequencies in the radial and angular motion.
Figure \ref{fig2} 
shows these POs in the RPL potential (\ref{potenra}) in the 
$(\tau,\alpha)$ plane, where
$\tau^{}_{\rm PO}=\tau(\alpha,L_{\rm PO})$ is the scaled period
for the PO specified by $(n_r,n_\theta)$, 
which satisfies the resonance
condition (\ref{rescond}).
As shown in this figure, the
polygonlike orbit $M(n_r,n_\theta)$ continues to exist 
after its
emergence at the bifurcation point $\alpha=\alpha_{\rm bif}$
from the parent circular 
orbit $MC$ ($M$th repetition of the primitive circle orbit $C$).
The exceptions are the diameter orbits $M$(2,1)
that exist for all values of $\alpha$
and form families with a higher degeneracy at the
HO symmetry-breaking point $\alpha=2$.

\subsection{Three-parameter PO families}
\l{sec5B:PISPM}

For the contribution of the three-parameter  ($\mathcal{K}=3$) 
families into the density 
shell correction $\delta g(E)~$ [Eq.~(\ref{pstrace})], 
one obtains \cite{MKApan16}
\bea\l{deltag3isp}
&&\delta g^{}_{{\rm scl}, P}(E) =  \Re \sum_{M P}
A_{M P}^{}(E) \nonumber\\ 
&&\times\exp\left[\frac{i}{\hbar}\,
S_{M P}(E)-i \frac{\pi}{2} \mu^{}_{M P}
-i\phi^{}_{\mathcal{D}}\right]\;.
\eea
The sum is taken over 
families of the three-parameter ($\mathcal{K}=3$)
polygonlike orbits $MP$ ($P$ stands for polygonlike, and $M$  
for the repetition number $M=1,2,...$).
In Eq.~(\ref{deltag3isp}), 
$S_{M P}(E)$ is the action along the PO,
\begin{equation}
S_{M P}(E) = 2\pi M \left[n_r\, 
I_r(E,L^*)+n_\theta \,L^*\right]\;,
\label{actionpo}
\end{equation}
where $I_r$ is the radial action variable
in the spherical 
phase-space coordinates
[Eq.~(\ref{ensurf})].
The angular momentum $L^\ast$ is given by the classical  
value for a particle motion along the $P$ orbit $L^\ast=L_{P}$ 
in an azimuthal plane.  
The numbers  $n_r$ and $n_\theta$ specify
the orbit $P$ with $n_r >2\,n_\theta$.  For the amplitude 
$A_{M P}^{}$ [Eq.~(\ref{deltag3isp})],
within the ISPM2, one finds
\be\l{amp3isp}
A_{M P}^{}=
\frac{L_{P} \,T_{P}}{\pi\hbar^{5/2} \sqrt{M n_r\,
K_{\rm P}}}\, \mbox{erf}\left(\mathcal{Z}_{M P}^{+},\mathcal{
  Z}_{M P}^{-}\right)
\,{\rm e}^{i \pi/4}\;,
\ee
where
$K_{P}$ represents the curvature (\ref{curvature}) at $L=L^\ast$  and
$T_{P}=T_{n_r,n_\theta}$ is the period of the
primitive ($M=1$) polygonlike orbit, $P(n_\theta,n_r)$.
For a three-parameter family
at the stationary point $L=L^*$ determined by 
the PO (stationary-phase) equation, one has 
\be\l{period3}
T_{P}=\frac{2\pi n_r}{\omega_r}=\frac{2\pi n_\theta}{\omega_\theta}\;.
\ee
The function $\mbox{erf}(\mathcal{Z}_{M P}^{+},\mathcal{
  Z}_{M P}^{-})$ in Eq.~(\ref{amp3isp}) is
given by the generalized error
function  (\ref{ispmlimits}).
Its complex arguments $\mathcal{Z}_{M P}^{\pm}$
are expressed in terms of the curvature $K_{P}$,
\bea\l{argerrorpar}
\mathcal{Z}_{M P}^{-} &=& \sqrt{-i \pi\, M\,n_r\, K_{P}/\hbar}\,\, 
\left(L_{-}-L_{P}\right)\;, \nonumber\\  
{\cal Z}_{M P}^{+} &=& \sqrt{-i\pi \,M\,n_r\, K_{P}/\hbar}
 \,\left(L_{+}-L_{P}\right)\;. 
\eea
We used here the ISPM2 
for the finite integration limits within the tori, i.e., between the
minimal $L_{-}=0$ and maximal $L_{+}=L_{C}$ values
of the angular-momentum integration variable for the ${\cal K}=3$ family 
contribution. The phase factor 
$\mu^{}_{M P}$ in Eq.~(\ref{deltag3isp}) 
is 
the Maslov index (see Sec.~\ref{sec3:cat3}) as
in the asymptotic Berry-Tabor trace formula. The amplitude (\ref{amp3isp}) 
obtained for the formula (\ref{deltag3isp}) is
regular at the bifurcation
points where the stationary point is located at the
end point $L=L^*=L_{C}$ of the action ($L$) part of a
torus.
 
For the case of the power parameter $\alpha$ sufficiently
far from the bifurcation points, one arrives at the SSPM limit
of the trace formula (\ref{deltag3isp})
with the amplitude $A_{M P}$, identical to the 
Berry-Tabor trace formula 
\cite{BT76} 
\begin{equation}
A_{M P}^{(\rm SSP)}= 
\frac{2\,L_{P} \,T_{P}}{\pi\hbar^{5/2} \sqrt{M n_r\,
K_{P}}}\, e^{i\pi/4}.
\label{amp3ssp}
\end{equation}
According to Sec.~\ref{sec3:cat3}, the Maslov index   
$\mu^{}_{M P}$ of Eq.~(\ref{deltag3isp}) 
is determined  by  
the number of
turning and caustic points within the FMCT 
(see Refs.~\cite{FEDoryuk_pr,MASLOV,FEDoryuk_book1,MAFptp06}),
\be\l{masl3Dra}
\mu^{}_{M P}
=3\,M\,n_r + 4\,M\,n_\theta\;, \qquad  
\phi^{}_{\mathcal{D}}=-\pi/2\;.
\ee
The total Maslov index $\mu_{M P}^{\rm (tot)}$  is defined as 
the sum of  
this asymptotic part (\ref{masl3Dra})
and the argument of the complex density amplitude (\ref{amp3isp})
\cite{migdalrev,MAFMptp02,MAFptp06,MFAMMSBptp99,MKApan16}. The total
index  $\mu_{M P}^{\rm (tot)}$
behaves as a smooth function of the energy $E$ and
the power parameter $\alpha$ through the bifurcation point.

\subsection{ Two-parameter circle families in the ISPM2}
\l{sec5C}

For contributions of the circular PO families to the trace formula
(\ref{pstrace}), one obtains
\bea\l{dengenC}
&&\delta g_{{\rm scl},C}^{}(E)=
\Re\sum_{M}\,A_{M C}^{} \nonumber\\ 
&\times&\exp\left[\frac{i}{\hbar}
\,S_{M C}\left(E\right) - \frac{i \pi}{2}
\mu^{}_{M C}-i\phi^{}_{\mathcal{D}}\right]\;. 
\eea
The sum is taken over  the repetition number for the circle 
PO and $M=1,2,...\;$. Here
$S_{M C}\left(E\right)$ is the action along the orbit 
$MC$, 
\be\l{actionC}
S_{M C}(E)=M \int_{0}^{2\pi} L \,\d\,\theta = 2\pi M\,L_{C}\;,
\ee
where $L_{C}$ is the angular momentum of the particle moving
along the orbit $C$.
For amplitudes of the $MC$-orbit
contributions, one obtains 
\bea\l{amp2ispC}
&&\!A_{M C}^{}(E)
= \frac{i L_{C}\,T_{C}}{ \pi\,\hbar^2\sqrt{F_{M C}}}\nonumber\\
&\times& \mbox{erf}\left(\mathcal{Z}_{p\;{M C}}^{+}\right)\,
 \mbox{erf}\left(\mathcal{Z}_{r\;{M C}}^{-},
\mathcal{Z}_{r\;{M C }}^{+}\right)\;,
\eea  
where $T_{C}=2 \pi/\omega^{}_{C}$ 
is the period of a particle motion along
the primitive ($M=1$) orbit C.  Here,  
$\omega^{}_{C}$ is the azimuthal frequency,
\be\l{omtcra}
\omega^{}_{C}
=\omega_\theta(L=L_{C})=L_{C}/(m\,r_{C}^2)\;, 
\ee
$r^{}_{C}$ the radius of the $C$ orbit, and 
$L_{C}$ is the angular momentum for a particle motion along the $C$ PO
\cite{KAprc12,KAijmpe04,MVApre13,MKApan16},
\be\l{rcLcd2Fcra}
\hspace{-0.5cm}r^{}_{C}\!=
\!R_0\left(\frac{2E}{(2\!+\!\alpha) V_0}\right)^{1/\alpha},\,\,\,\,\, 
L_{C}\!=\!p(r^{}_{C})r^{}_{C}\;.  
\ee
In Eq.~(\ref{amp2ispC}), $F_{M C}$ is the
stability factor (the trace of the monodromy matrix)
\be\l{fgutz}
F_{M C}= 4 \sin^2\left(\pi M\sqrt{\alpha+2}\right)
\ee
and  $\mathcal{J}_{M C}^{(p)}$ is the Jacobian
\be\l{jacpCres}
\mathcal{J}_{M C}^{(p)}= 2 \pi (\alpha+2)\;
M\; K_{C}\;r_{C}^{2}\;,
\ee
where $K_{C}$ is the curvature
for $C$ orbits \cite{MVApre13},
\be\l{curvraC}
K_{C} = -\frac{(\alpha+1)(\alpha-2)}{12 \,(\sqrt{\alpha+2})^3\,
  L_{C}}.  
\ee
The finite limits in the error functions 
of Eq.~(\ref{amp2ispC}), ${\mathcal Z}_{p\;M C}^{\pm}$ and 
${\mathcal Z}_{r\;M C}^{\pm}$,
are given by 
\bea\l{ZplimCfin}
\!{\mathcal Z}_{p\;M C}^{-}\!&=&\!0,\,\,\,\,\, 
{\mathcal Z}_{p\;M C}^{+}\!=\!L_{C} \sqrt{-\frac{i \pi}{\hbar}
(\alpha+2) M K_{C}}\;,\nonumber\\
{\mathcal Z}_{r\;M C}^{-}&=&
\left(\frac{r_{\rm min}}{r^{}_{\rm C}}-1\right)\;\sqrt{\frac{i F_{M C}}{4 \pi 
(\alpha+2)\,\hbar M K_{C}}}\;,\nonumber\\
{\mathcal Z}_{r\;M C}^{+}&=&
\left(\frac{r_{\rm max}}{r^{}_{C}}-1\right) \sqrt{\frac{i F_{M C}}{4 \pi 
(\alpha+2)\,\hbar M K_{C}}}\,, 
\eea
 where $r_{\rm min}$ and $r_{\rm max}$ are the radial turning
points specified below. The interval between them covers the CT manifold 
including  
the stationary point $r^{}_{C}$.

Asymptotically far 
from the bifurcations 
(also far from the symmetry 
breaking point $\alpha=2$), 
the amplitude (\ref{amp2ispC}) approaches the SSPM limit
\begin{equation}
A_{M C}^{}(E) \rightarrow 
A_{M C}^{\rm SSP}(E)
= 
\frac{2 L_{C}\,T_{C}}{\pi \hbar^2\;\sqrt{F_{M C}}}\;.
\label{amp2sspC}
\end{equation} 
In these SSPM derivations, the radial integration limits 
$r_{\rm min}$ and 
$r_{\rm max}$  
turn into the asymptotic values
\be\l{rmimax}
r_{\rm min}=0, \qquad 
r_{\rm max}=  
R_0\mathcal{E}^{2/(2+\alpha)}\;.
\ee
They are given by 
the two real solutions of Eq.~(\ref{pr2eq}) at $L=0$.  
The upper limit $r_{\rm max}$ was
extended to infinity 
in the derivation of Eq.~(\ref{amp2sspC})
because the stationary point $r^{}_{\rm C}$ is far away from 
both integration boundaries $r_{\rm min}$ and  $r_{\rm max}$. 

The Maslov index $\mu^{}_{M C}$ 
in Eq.~(\ref{dengenC}) is given by
\begin{equation}
\mu^{}_{M C} 
=2\,M\;,\qquad 
\phi^{}_{\mathcal{D}}=\pi/2\;.
\label{masl2}
\end{equation}
For the calculation of this asymptotic 
Maslov  index $\mu^{}_{M C}$ through
the turning and caustic points [see the trace formula 
(\ref{dengenC}) with the ISPM2 
(\ref{amp2ispC}) and the SSPM  
(\ref{amp2sspC}) amplitude], one can use the FMCT 
(Sec.~\ref{sec3:cat3}). The total Maslov index  $\mu_{M C}^{\rm (tot)}$
can be introduced as above (see 
Refs.~\cite{migdalrev,MAFMptp02,MAFptp06,MFAMMSBptp99,MKApan16}).

Taking the opposite limit 
$\alpha \rightarrow \alpha_{\rm bif}=n_r^2/n_\theta^2-2$  to the bifurcations 
where $F_{M C} \rightarrow 0$,
but far away from 
the HO limit $\alpha=2$, 
one finds that the argument of the second  error function
in Eq.~(\ref{amp2ispC}), coming from the radial-coordinate integration, 
tends to zero proportional to $\sqrt{|F_{M C}|}$  
[Eq.~(\ref{ZplimCfin})].
 Thus, as in a general case (Sec.~\ref{sec4:bifurcations}), the singular  
stability 
factor $\sqrt{F_{M C}}$ of 
the denominator in Eq.~(\ref{amp2ispC}) is exactly canceled by 
the same from the numerator.  At the bifurcation $F_{M C}\to
0$, one obtains
\bea\l{amp2ispCbif}
A_{M C}^{\rm ISP}(E)
&\rightarrow& \frac{2L_{C}}{\hbar^{3/2}\; \omega^{}_{C}\;
\sqrt{i\pi (\alpha+2)\,M\,K_{C} r^2_{C}}}\nonumber\\
&&\times\mbox{erf}\left({\cal Z}_{p\;{M C}}^{(+)}\right)
\left(r^{}_{\rm min}-r^{}_{\rm max}\right)\;.
\eea
Therefore,
Eq.~(\ref{amp2ispC}) gives the finite result through the bifurcation.
Taking the reduction  
of the end-point manifold for the parent-orbit
term, according to Eq.~(\ref{bifboundlim}),
\be\l{rboundbif}
r^{}_{\rm min} \rightarrow r^{}_{\rm max}\rightarrow r^{}_{\rm C}\;,
\ee
the amplitude (\ref{amp2ispC}) vanishes
at the bifurcation point \cite{MKApan16}  
[see Eq.~(\ref{amp2ispCbif})].
This is in line with general arguments for the bifurcation
limit [see Sec.~\ref{sec4:bifurcations} around Eq.~(\ref{bifboundlim})].

As shown in the Appendix,  
following the FMCT 
(Sec.~\ref{sec3:cat3}) one can derive
the ISPM3 expression for the oscillating level density.
One may note that the parameter $w$ given by Eq.~(\ref{wC})
can be considered as the
dimensionless semiclassical measure of the distance from a bifurcation.
Similarly, as for the caustic catastrophe
points, in the case of the application of the FMCT 
(Sec.~\ref{sec3:cat3}) to the bifurcation of the circular
orbits (see also Sec.~\ref{sec4:ispm}), one can use 
simply the ISPM2 (Sec.~\ref{sec5C})
as the simplest approximation 
near the bifurcation, i.e., $w \siml 1 $. 
However, working out properly 
the transition itself from the asymptotic
SSPM radial-integration limits $r^{}_{\pm}$  to 
the same value $r^{}_{C}$ at the bifurcation
within its close vicinity ($w \siml 1 $) 
is left for future work. 

\subsection{Two-parameter diameter families}
\l{sec5D:DSSPM}

For the diameter-orbit ($\mathcal{K}=2$) family contribution
to the trace formula (\ref{pstrace}) for the RPL potential, the ISPM is needed
only near the symmetry-breaking at $\alpha=2$ of
the harmonic-oscillator limit \cite{MKApan16}.   
For our purpose,  
we simply use the SSPM approximation for the diameter families,
valid at the values of the power parameter $\alpha$ far from
the symmetry-breaking limit,
\bea\l{dengenD}
&&\delta g_{{\rm scl}\;D}^{}(E)=
\Re\sum^{}_{M} A_{M D}^{}\nonumber\\
&&\qquad\times
\exp\left[\frac{i}{\hbar}
S_{M D}\left(E\right) - \frac{i\pi}{2}
\mu^{}_{M D}-i\phi_{\mathcal{D}}\right],
\eea
where
\be\l{amp2sspD}
A_{M D}^{} 
=\frac{1}{i \pi M K_{D}\omega_r \hbar^2}\;.
\ee
The frequency $\omega_r$ is expressed through the radial period,  
\be\l{periodfreq} 
T_r=\frac{2 \pi}{\omega_r}=\int_0^{r_{\rm max}}\frac{2m\d r}{p(r)}
=\sqrt{\frac{2m\pi}{E}}
\frac{r_{\rm \max} \Gamma\left(1+1/\alpha\right)}{ 
\Gamma\left(1/2 + 1/\alpha\right)}\;,
\ee
where $\Gamma(x)$ is the Gamma function.
In Eq.~(\ref{amp2sspD}), $K_{D}$ is the diameter curvature
\cite{MVApre13}
 \be\l{curvraD}
K_{D}=
\frac{\Gamma\left(1-1/\alpha\right)}{\Gamma\left(1/2-1/\alpha\right)\epsi 
\;\sqrt{2 \pi m R_0^2 V_0}}\;.
\ee
For the Maslov index [Eq.~(\ref{dengenD})], one obtains
\be\l{maslD}
\mu^{}_{M D} 
=2M\;,
\qquad \phi^{}_{\mathcal{D}}=-\pi/2\;.
\ee

\subsection{TOTAL TRACE FORMULAS FOR THE SPHERICAL RPL POTENTIAL}
\l{sec5E:total}

The total ISPM trace formula for the RPL potential
is the sum of the contribution of
the $\mathcal{K}=3$ polygonlike ($P$) families $\delta g_{P}^{}(E)$
[Eqs.~(\ref{deltag3isp}) with (\ref{amp3isp})],
the $\mathcal{K}=2$ circular ($C$) families
$\delta g_{C}^{}(E)$ [Eqs.~(\ref{dengenC}) and (\ref{amp2ispC})
for the ISPM2  and Eqs.~(\ref{denC3}) and (\ref{ampc3}) for the ISPM3], 
and the $\mathcal{K}=2$ diameter ($D$) families
$\delta g_{D}^{}(E)$
[Eqs.~(\ref{dengenD}) and (\ref{amp2sspD})], 
\be\l{deltadenstotprlp}
\delta g_{\rm scl}(E)
=\delta g_{{\rm scl},P}^{}(E) + \delta g_{{\rm scl},C}^{}(E) +
\delta g_{{\rm scl},D}^{}(E)\;.
\ee
This trace formula has the correct finite asymptotic limits to
the SSPM: The Berry-Tabor result [Eqs.~(\ref{deltag3isp}) 
and (\ref{amp3ssp})] 
for the $P$ orbits ($\mathcal{K}=3$)
and for $C$ orbits ($\mathcal{K}=2$)
[Eqs.~(\ref{dengenC}) and (\ref{amp2sspC})]; see the same 
for $D$ orbits [Eqs.~(\ref{dengenD}) and (\ref{amp2sspD})].
Transforming the variable from the ordinary energy $E$ to the
scaled energy $\epsi$ (\ref{eq:scaledentau1}), one obtains the trace
formula for the scaled-energy level density,
\be\l{scldenstyra1}
\delta \mathcal{G}(\epsi)=\delta g(E)\frac{dE}{d\epsi}
=\sum_{\rm PO}\delta\mathcal{G}_{\rm PO}(\epsi),
\ee
with
\begin{align}
\delta\mathcal{G}_{\rm PO}(\epsi)
=&\Re\left[\mathcal{A}_{\rm PO}(\epsi)\right. \nonumber \\ &\times
\left.\exp\left(i\epsi\tau_{\rm PO}-\frac{i\pi}{2}\mu_{\rm PO}
-i\phi_{\mathcal{D}}\right)\right]. \l{scldenspo}
\end{align}
The Fourier transform of this scaled-energy level density,
truncated by the
Gaussian with the cutoff $\gamma$, is expressed as
\begin{align}
F_\gamma(\tau)&=\int\mathcal{G}(\epsi)e^{i\epsi\tau}e^{-(\epsi/\gamma)^2}
d\epsi
\nonumber\\
&=\sum_{\rm PO}\tilde{\mathcal{A}}_{\rm PO}(\tau)
e^{-\gamma^2(\tau-\tau_{\rm PO})^2/4}.
\l{ftlscl}
\end{align}
This gives a function with successive peaks at the scaled periods of the
classical POs $\tau=\tau_{\rm PO}$ with the height $|\tilde{\mathcal{A}}_{\rm
PO}|$, which is proportional to the amplitude $\mathcal{A}_{\rm PO}$ of the
contribution of the orbit PO to the semiclassical level density.
Evaluating the same Fourier
transform by the exact quantum level density, one has
\bea\l{fourierpower1}
F(\tau)&=&\int \left[\sum_i\delta(\epsi-\epsi_{i})\right]
e^{i\epsi \tau}e^{-(\epsi/\gamma)^2}\d\epsi \nonumber\\ 
&=&\sum_i e^{i\epsi_{i}\tau}e^{-(\epsi_i/\gamma)^2},
\eea 
with
\be\l{epsi}
\epsi_{i}=\left(E_i/V_0\right)^{1/2+1/\alpha}\;.
\ee
Thus, one can extract the contribution of
classical periodic orbits to the level density from
the Fourier
transform of the quantum level density.
In what follows we consider the classical-quantum correspondence
using this Fourier transformation technique, in addition to the
direct comparison of quantum and semiclassical level densities.

\section{Comparison with quantum results}
\label{sec6:comp-quant}

\begin{figure}[bt]
\begin{center}
\includegraphics[width=\linewidth,clip=true]{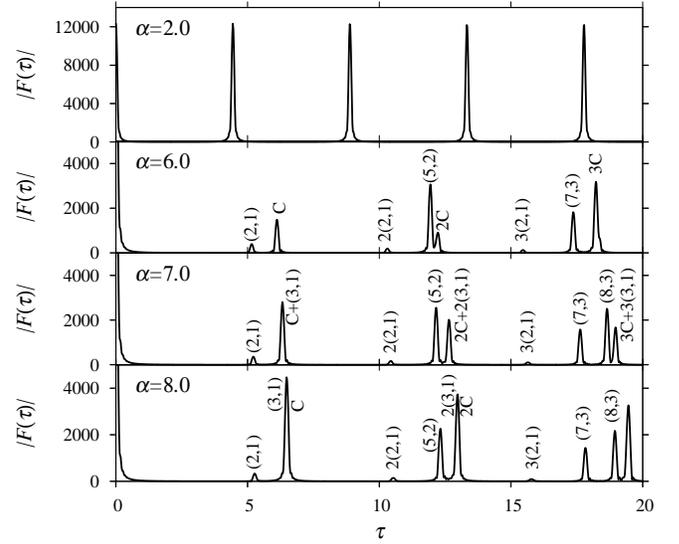}
\end{center}
\caption{Moduli of the Fourier transform $|F(\tau)|$ of 
the quantum scaled-energy level density [Eq.~(\ref{fourierpower1})] 
plotted for several values of $\alpha$.
}
\label{fig3}
\end{figure}

Figure~\ref{fig3} 
shows the Fourier transform of the 
quantum-mechanical level density $\mathcal{G}(\epsi)$ for the RPL
potential [see Eq.~(\ref{fourierpower1})].
At the HO limit $\alpha=2$, all the classical orbits are periodic
and form the four-parameter family for a given energy.  The Fourier
transform exhibits the
equidistant identical peaks at $\tau_n=\sqrt{2}\pi n$,
corresponding to the $n$th repetitions of the primitive PO family.
With increasing $\alpha$,  
each peak is split into two peaks corresponding to the diameter ($D$) and
circle ($C$) orbits and
the amplitudes of the oscillating level density for these  
orbits are decreased.
However, one finds a growth of the peak at $\tau\sim 6.2$
corresponding to the $C$ orbit
around the bifurcation point 
$\alpha^{}_{\rm bif}=7.0$.
Note that,  
approaching the bifurcation, 
the contribution of the
$C$ orbit is strongly enhanced until it forms a local family of POs
with a higher degeneracy at the bifurcation point. 
From this point a trianglelike 
$P$(3,1) family bifurcates.
This family has  
high degeneracy $\mathcal{K}=3$.
It remains important, also  
for larger $\alpha >\alpha^{}_{\rm bif}$.
The above enhancement in the Fourier peaks $F(\tau)$
is directly associated with the oscillating ISPM
level-density amplitude $A^{}_{\rm PO}$ of
the bifurcating PO family having a high degeneracy.
This family  
is a major term in the $\hbar$ expansion
in the comparison with the SSPM asymptotics 
(see Sec.~\ref{sec4:bifurcations} and Ref.~\cite{MKApan16}). 
The Fourier peak at $\tau \sim 6.2$ in Fig.~\ref{fig3}
shows the enhancement of the amplitude of the newborn 
$P$(3,1) family contribution including $C$(1,1) orbits as the end points
(see the Introduction and Sec.~\ref{sec4:bifurcations}).

\begin{figure}[tb]
\begin{center}
\includegraphics[width=\linewidth,clip=true]{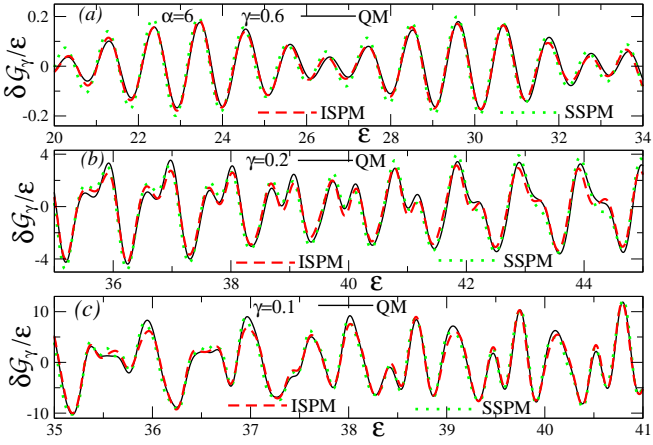}
\end{center}
\caption{Comparison of the  quantum-mechanical (QM, black solid line) 
and semiclassical (ISPM, red dashed line, and SSPM, green dotted line)
shell-correction scaled-energy
level density $\delta \mathcal{G}_\gamma(\epsi)$ [(\ref{scdentot}) and
(\ref{qmdentot})], divided by $\epsi$,
as a function of the scaled-energy $\epsi$ for $\alpha=6.0$ 
and averaging the
parameters  (a) $\gamma=0.6$,  
(b) $\gamma=0.2$,  
and (c) $\gamma=0.1$. 
}
\label{fig4}
\end{figure}

\begin{figure}[tb]
\begin{center}
\includegraphics[width=\linewidth,clip=true]{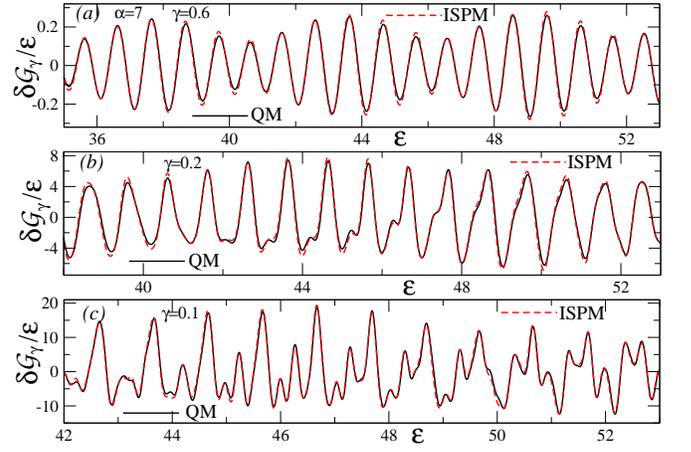}
\end{center}
\caption{Same QM and ISPM as in Fig.~\ref{fig4} but  
for $\alpha=7.0$. 
}
\label{fig5}
\end{figure}

\begin{figure}[tb]
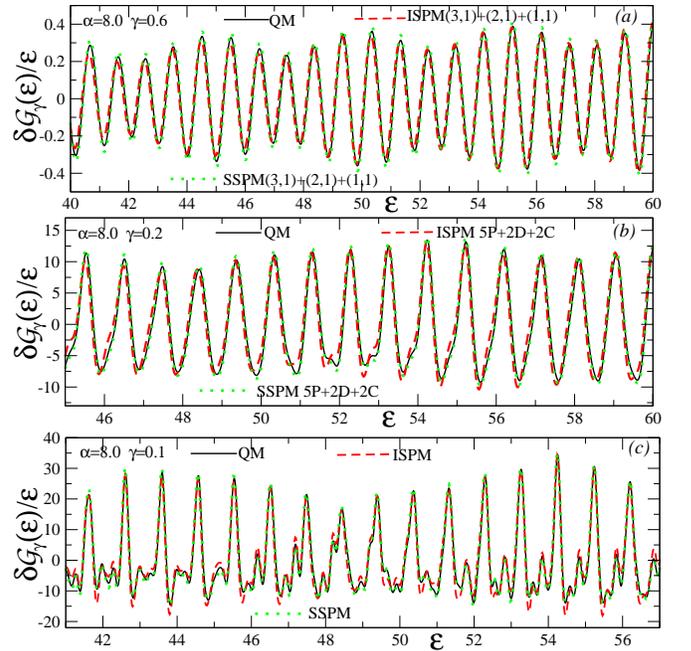

\begin{center}
\includegraphics[width=\linewidth,clip=true]{fig6a_dg3D_a800_g06_1P1D1C.eps}
\includegraphics[width=\linewidth,clip=true]{fig6b_dg3D_a800_g02_5P2D2C.eps}
\includegraphics[width=\linewidth,clip=true]{fig6c_dg3D_a800_g01_FULLPDC.eps}
\end{center}
\caption{Same as in Figs. \ref{fig4} and \ref{fig5} but for 
$\alpha=8.0$: 
(a) major shell structure at $\gamma=0.6$, 
(b) fine shell structure at $\gamma=0.2$, 
and (c) fine shell structure at $\gamma=0.1$. The shortest
POs and their numbers taken into account in the calculations are shown
in the legends.
}
\label{fig6}
\end{figure}

As the significance of bifurcations is confirmed through the
Fourier analysis [Eqs.~(\ref{ftlscl}) and (\ref{fourierpower1})], 
let us now investigate the oscillating part of the
scaled-energy level densities averaged
with the Gaussian averaging parameter $\gamma$,
\be\l{avdentot}
\delta\mathcal{G}_\gamma(\epsi)
=\int \exp\left[-\left(\frac{\epsi-\epsi'}{\gamma}\right)^2\right]~
 \delta\mathcal{G}(\epsi')d\epsi'.
\ee
The semiclassical shell-correction density
$\delta\mathcal{G}_{\gamma,\;{\rm scl}}(\epsi)$ is given by 
\be\l{scdentot}
\delta\mathcal{G}_{\gamma,\;{\rm scl}}(\epsi)
=\sum_{\rm PO}\delta\mathcal{G}_{\rm PO}(\epsi)
 \exp(-\tau_{\rm PO}^2\gamma^2/4)\;.
\ee
[see Eq.~(\ref{scldenspo}) for $\delta \mathcal{G}_{\rm PO}$].
For the quantum density, one has
\be\l{dGQM}
\delta \mathcal{G}_{\gamma,{\rm QM}}=\mathcal{G}_{\gamma,{\rm QM}}-
\widetilde{\mathcal{G}}_{\rm QM}\;,
\ee
where
\be\l{qmdentot}
\mathcal{G}_{\gamma,\;{\rm QM}}(\epsi)
=\sum_i \exp\left[-\left(\frac{\epsi-\epsi_i}{\gamma}\right)^2\right]\;.
\ee
The smooth level density $\widetilde{\mathcal{G}}_{\rm QM}$ is calculated
for the scaled spectrum $\epsi_i$.
For these calculations we employed 
the standard Strutinsky averaging
(over the scaled energy $\epsi$), finding a good plateau\footnote{
It may worth pointing out that the 
quality of the plateau in the SCM calculations
of the level density  
is much better  
when using the scaled-energy variable $\epsi$ rather 
than the energy $E$
itself \cite{MVApre13}.}
around the 
Gaussian averaging width ${\widetilde \gamma}=2-3$ and curvature-correction 
degree ${\cal M}=6$.

Figures \ref{fig4}--\ref{fig6}  
show good agreement of the coarse-grained ($\gamma=0.6$)
and fine-resolved ($\gamma=0.1-0.2$) semiclassical
and quantum results for
$\delta \mathcal{G}^{}_\gamma(\epsi)$ 
(divided by $\epsi$) as functions 
of the scaled energy $\epsi$ at  
$\alpha=6.0, 7.0$, and $8.0$.
At $\alpha=6$, as well as $\alpha=2$ and $4$,
the analytic expressions
for all of the classical PO characteristics
are available, and can be used to check the precision of
the numerical calculations \cite{MVApre13,MKApan16}.
The values $\alpha=6$ (Fig.~\ref{fig4}) and $8$ (Fig.~\ref{fig6})
are taken as examples that 
are sufficiently far 
from the bifurcation point $\alpha=\alpha_{\rm bif}=7$
(Fig.~\ref{fig5}).
The ISPM results at these values of $\alpha$
show good convergence to
the SSPM results.
The $C$ and $D$
POs with the shortest (scaled) periods $\tau$ are dominating
the PO sum
at a large averaging parameter $\gamma=0.6$ (the 
coarse-grained or 
major-shell structure). Many 
more families with a relatively
long period $\tau$ at $\gamma=0.1-0.2$ (the fine-resolved shell structure)
become significant 
in comparison with the quantum results 
\cite{MKApan16}.

Figure~\ref{fig5} shows the results for the bifurcation point
$\alpha=7$ where
the trianglelike (3,1) PO family emerges from the parent $C$(1,1)
family in a typical bifurcation scenario.  One also finds
good agreement of the ISPM with quantum results here.
As the SSPM approximation fails at 
the bifurcation, 
it is not presented
in Fig.~\ref{fig5}.
As discussed above, the SSPM
at the bifurcation yields a sharp
discontinuity of the $P$(3,1) amplitude and a divergent behavior of the 
$C$(1,1), in contrast to the continuous ISPM components.
Our results in Fig.~\ref{fig5} demonstrate
that the ISPM
successfully solves these
catastrophe problems of the SSPM for all
averaging parameters $\gamma$.
In contrast to the results
shown in Figs.~\ref{fig4} and \ref{fig6},
the coarse-grained ($\gamma=0.6$) density oscillations at the
bifurcation point $\alpha=7$ (Fig.~\ref{fig5})
do not contain any contributions
from the $C$(1,1)
end-point term but instead the $P$(3,1) term becomes 
dominant.

Note that  
many more families with relatively
long periods $\tau$ become necessary to account for the
fine-resolved shell structures ($\gamma=0.1-0.2$) \cite{MKApan16}. 
For the exemplary bifurcation $\alpha=7.0$, at smaller
averaging parameters ($\gamma \siml 0.2$) the dominating orbits 
become the bifurcating newborn 
$P$(3,1) of the highest degeneracy $\mathcal{K}=3$ 
along with the leading $P$(5,2), $P$(7,3), and $P$(8,3) POs 
which are born at smaller $\alpha$ (see Fig.~\ref{fig2}).
They include the parent $C$-orbit end-point manifolds.
As also shown in
the quantum Fourier transforms in Fig.~\ref{fig3},
these POs yield larger 
contributions at the bifurcation values of $\alpha$ and are even more enhanced 
on their right in a wide region of $\alpha$.

\section{Conclusions}
\l{sec6:concl}

The Fedoriuk-Maslov catastrophe theory is extended to simple
bifurcation problems in the POT.
Within the extended FMCT,
we overcome the divergence and discontinuity
of semiclassical amplitudes of the 
standard stationary-phase method, in particular, in
the Berry-Tabor formula
near bifurcations. 
A fast convergence in the PO expansion of the averaged level density 
for a large Gaussian averaging parameter is shown too.
This allows one often to express significant 
features of the shell structure in terms of a 
few short periodic orbits.
We have formulated our ISPM trace formula for a simple
bifurcation scenario so that the parent orbits at the end points
have vanishing contributions
at the bifurcation point,  
which allows us to consider them 
everywhere separately from the term for a
newborn family of the periodic orbits.

The extended FMCT is used for derivations of
the trace formula in the case of
the three-dimensional spherical RPL potential
by employing the improved stationary phase method.
We presented a class of the radial power-law
potentials that, 
up to a constant,
provides a good approximation to the WS potential
in the spatial region where the particles are bound. 
The RPL potential is
capable of controlling surface diffuseness and  
contains
the popular harmonic-oscillator and cavity potentials in
the two limiting cases of the power parameter $\alpha$.
Its advantage is
the scaling invariance of
the classical
equations of motion.
This invariance makes the POT calculations and the Fourier analysis of the
level density very easy.
The contribution of the POs to the semiclassical level density
and shell energies is expressed analytically
(and even all the PO characteristics are given explicitly,
e.g., for $\alpha=6$)
in terms of the simple special functions.
The quantum Fourier spectra yield directly the amplitudes of the
quantum level density at the periods (actions) of the
corresponding classical POs.

We have derived the semiclassical trace formulas that   
are also valid in the bifurcation region and examined them at
the bifurcation catastrophe points and asymptotically far
from them in the spherical RPL potential model.
They are based on   
the SPM improved
to account for the effect of the bifurcations by using the extended
FMCT.
The ISPM overcomes the problems of singularities in the SSPM and
provides the generic trace formula that relates
the oscillating component of the level density for a 
quantum system to a sum over POs of the corresponding 
classical system.
We showed good convergence
of this improved trace formula to the simplest ISPM based on the second-order
expansion of the
classical action 
at several characteristic values of 
the power parameter $\alpha$ including the bifurcations
and asymptotically far from them. 

We obtained good agreement between the 
ISPM semiclassical and quantum results for the 
level-density shell corrections 
at different values of the power parameter $\alpha$,
both at the bifurcations and far from them.
Sufficiently far from the bifurcation of the leading
short POs with a maximal degeneracy,  
one finds also  
good convergence
of the ISPM trace formulas to the SSPM  approximation. 
We  emphasize the 
significant influence 
of the bifurcations of short POs on the main characteristics 
of oscillating components of the single-particle
level density for a fermionic 
quantum system.  They appear
in the significant fluctuations
of the energy spectrum (visualized 
by its Fourier transform), namely, the shell structure. 

In line with the general arguments of the extended FMCT, 
the stationary points forming the circular-orbit families
are located at the end point of the classically accessible region
and they coincide with
the newborn family of the polygonlike orbits at the bifurcation.
Taking into account the  
reduction of the end-point manifold in the
bifurcation limit, the parent $C$-family contribution is transformed
into the newborn $P$-family term that
presents now their common result.
Thus, one has the separate contributions of 
the parent $C$ and newborn $P$
orbits through the bifurcation scenario, but with no concern about  
double counting.

Future work  
should study in detail the 
transition
of the ISPM trace formula from the bifurcation points 
to its asymptotic SSPM region.
This will enable us to understand more properly the shape
dynamics of the finite fermion systems.  In particular, the improved
stationary phase method
can be applied to describe
the deformed shell structures
where bifurcations play  
an essential role
in formations of the
superdeformed minima along a
potential energy valley
\cite{migdalrev,MKApan16}. One
of the remarkable tasks might be to clarify, in terms of the
symmetry-breaking (restoration) and bifurcation phenomena, 
the reasons of the exotic
deformations such as  the octupole and tetrahedral ones
within the suggested  ISPM.  In this way, it would be worth extending
our present local bifurcation FMCT  to describe, e.g.,  a bridge
(non-local) bifurcation phenomenon 
found in a more realistic mean field in the fermionic systems (see also 
Refs.~\cite{KAprc12,AMprc14,KAps16,KAps17}).

Our semiclassical analysis may therefore lead 
to a deeper understanding of the shell effects in the
finite fermionic systems such as atomic nuclei, metallic clusters, 
trapped fermionic atoms, and 
semiconductor quantum dots
\cite{sclbook,ABjpa08,RBMpra96,BBCzp97,FKMSprb08,BO05}.
Their level densities, conductance, and magnetic 
susceptibilities
are significantly modified by shell effects. 
As a first step towards the collective
dynamics, the oscillating parts of the nuclear moment of inertia 
should be studied semiclassically in terms of POs taking into account the
bifurcations \cite{MSKBprc10,MGBpan14,GMBBprc16,MGBpan17}.

\section*{ACKNOWLEDGMENTS}

The authors gratefully acknowledge
M.\ Brack and K.\ Matsuyanagi
for fruitful collaborations and 
many useful discussions.  
One of us (A.G.M.) is
also very grateful for  
hospitality during his working visits to the 
Department of Physical Science and Engineering of the Nagoya Institute of 
Technology and for financial support from 
the Japanese Society of Promotion of Sciences 
through Grant No. S-14130.

\appendix

\section{The ISPM3 approximation}

Following the FMCT (Sec.~\ref{sec3:cat3}),
one can derive
the improved (ISPM3) contribution of the $C$ orbits
by taking
into account the third-order
terms of the action
expansion [see Eq.~(\ref{Sexp})] in the integration 
over the catastrophe variable $r$.
Within the ISPM3, one obtains 
\begin{align}\l{denC3}
&\delta g_{{\rm scl},C}^{(3)}(E)=
\Re\sum_{M}\,A_{M C}^{(3)} \nonumber\\ 
&\times\exp\left[\frac{i}{\hbar}
\,S_{M C}\left(E\right) - \frac{i \pi}{2}
\mu^{}_{M C}+\frac{2i}{3}w^{3/2}-i\phi^{}_{\mathcal{D}}\right].
\end{align}
The ISPM3 amplitudes $A_{M C}^{(3)}$ are given by
\begin{align}\l{ampc3}
A_{M C}^{(3)}
=& \frac{2\Lambda\; \sqrt{L_{C}}}{\hbar^{5/2}\omega^{}_{C}}\;
\frac{\mbox{erf}\left(\mathcal{Z}^{+}_{p,M C}\right)}{
\sqrt{2 \pi i (\alpha+2)M K_{C}}} \nonumber\\
&\times\left[\mbox{Ai}\left(-w,\mathcal{Z}_{M C}^{(-,3)},
\mathcal{Z}_{M C}^{(+,3)}\right) \right.\nonumber\\
&+\left. 
i\,\mbox{Gi}\left(-w,\mathcal{Z}_{M C}^{(-,3)},\mathcal{Z}_{M C}^{(+,3)}
\right)\right]\;,
\end{align} 
where
\be\l{wC}
w=w^{}_{M C}=\left[\frac{\kappa^{1/3}\;F_{M C}}{
4 \pi (\alpha+2) M K_{C} L_{C}(3a)^{2/3}}\right]^2\;,
\ee
with 
\be\l{akS}
a=\frac{r^{3}_{C}}{6L_{C}}\; \Phi_{CT}'''(r^{}_{C})\;,\quad
\quad 
\kappa=\frac{L_{C}}{\hbar}\;.
\ee
The parameter $\epsilon$
of Eq.~(\ref{catint}) [Eq.~(\ref{eps})] used in deriving the above
expressions 
is proportional
to the stability factor $F_{M C}$ [Eq.~(\ref{fgutz})],
\be\l{epsC}
\epsilon= 
\frac{F_{M C}}{4 \pi (\alpha+2) M 
K_{C}L_{C}}\;.
\ee
The incomplete Ai (Gi) integrals in Eq.~(\ref{ampc3})
are defined by Eq.~(\ref{airynoncomp}).
The integration limits of these 
functions are the same as those given by  Eq.~(\ref{zpm}),
\be\l{limrC3}
\mathcal{Z}_{M C}^{(\pm,3)}=
\frac{r^{}_{\pm}-r^{}_{C}}{r^{}_{C}\Lambda}
+\sigma \sqrt{w^{}_{M C}}\;, \quad
\sigma=\sgn(\epsilon)\;,
\ee 
where $r_\pm$
are the upper ($r^{}_{+}>r^{}_{C}$) and  
lower 
($r^{}_{-}<r^{}_{C}$) limits for the radial integration.
These integration limits are defined in Eq.~(\ref{rmimax}).
In the bifurcation limit $\alpha \rightarrow \alpha^{}_{\rm bif}$,
both terms in  
square brackets  
in Eq.~(\ref{ampc3}) [see also Eq.~(\ref{limrC3})]
go to zero  
for the same reason as in the ISPM2 case. 
The incomplete Ai and Gi functions of the 
integrand [Eq.~(\ref{airynoncomp})] have no singularities in the 
bifurcation limit $w \rightarrow 0$. 
In addition, 
the radial integration limits $r_{\rm min}$ and  $r_{\rm max}$ 
approach the stationary point $r^*=r^{}_{C}$,
which is the $C$-orbit radius [Eq.~(\ref{rboundbif})].
This ensures the disappearance of the end-point manifold in this limit
and therefore,
in line with the general arguments of
Sec.~\ref{sec4:bifurcations} (see also Sec.~\ref{sec5C}),  
one finds the zero contribution of 
the circular orbit term exactly at this bifurcation.
In turn,  
the contribution of the circular orbit is included 
in the newborn $P$ orbit term.

In the opposite limit,
sufficiently far from the bifurcation points, where the
stability factor
$F_{M C}(\alpha)$ takes a finite nonzero value,
the second term in Eq.~(\ref{limrC3}) changes
with increasing $\epsilon$
much faster (proportional to $\kappa^{1/3}\epsilon$)
than the first component (proportional to $\kappa^{1/3}$).
Thus, one has
$|\mathcal{Z}_{M C}^{(\pm,3)}|\gg |z^{\ast,\pm}_{M C}|$ 
in this limit.
The ISPM3 expression 
[Eq.~(\ref{denC3}) with Eq.~(\ref{ampc3})] for the oscillating level
density suggests that the parameter $w$
($w \propto F_{M C}^2 \propto \epsilon^2$)
can be considered as  
a dimensionless measure of the distance from the bifurcation
[see Eqs.~(\ref{wC}) and (\ref{epsC})].
For a large distance from the 
bifurcation  $w \gg 1$ [Eq.~(\ref{wC})],
one can extend the radial integration limits
as $r_-=0$ and $r_+\to\infty$.
In this limit, the incomplete Airy and Gairy
functions (\ref{airynoncomp}) can be approximated by
the complete ones 
(\ref{airystandard}). Since the argument $w$ of these standard
functions [Eq.~(\ref{wC})] at a finite stability factor 
$F_{M C}(\alpha)$ becomes large in the semiclassical 
limit $\kappa \gg 1$, 
one can use their asymptotic expressions (\ref{AiGiAsympt}).
Thus,  we arrive at the same SSPM result 
[Eqs.~(\ref{dengenC}) and (\ref{amp2sspC})] for the $C$ family contributions
as obtained from 
the ISPM2 $C$-trace formula  [see Eq.~(\ref{ampc3})
for its amplitude]. 

For the simplest catastrophe problem, the ISPM3 might become important
when the PO is distant from the bifurcation points to some extent but not
asymptotically far from them.  It is also necessary for the higher-order
catastrophe problem, which is  
not found in the RPL model discussed in this paper.
The definition of the end-point manifold might also be 
affected by the
consideration of  
higher expansion terms.  This is also a problem
to be solved in the future in order to describe the transition from
a bifurcation vicinity to the asymptotic  region.


\end{document}